\input amstex
\documentstyle {amsppt}

\pagewidth{32pc} 
\pageheight{45pc} 
\mag=1200
\baselineskip=15 pt

\hfuzz=5pt 
\topmatter
\NoRunningHeads 
\title finite localities I 
\endtitle
\author Andrew Chermak
\endauthor
\affil Kansas State University 
\endaffil
\address Manhattan Kansas
\endaddress
\email chermak\@math.ksu.edu
\endemail
\date
November 2021  
\enddate 

\endtopmatter

\define\w{\widetilde}
\define\wh{\widehat}
\redefine\norm{\trianglelefteq}

\redefine\bar{\overline}

\redefine\maps{\mapsto}
\redefine\i{^{-1}}

\redefine\l{\lambda}
\redefine\s{\sigma}
\redefine\a{\alpha}
\redefine\b{\beta}

\redefine\g{\gamma}

\redefine\r{\rho}

\redefine\G{\Gamma}

\redefine\S{\Sigma}
\redefine\L{\Lambda}

\redefine\<{\langle}
\redefine\>{\rangle}

\redefine\ca{\Cal}

\redefine\D{\Delta}

\redefine\sub{\subseteq} 

\redefine\nsub{\nsubseteq}

\redefine\nset{\emptyset} 

\redefine\1{\bold 1} 

\redefine\up{\uparrow}

\document

\vskip .1in 
\noindent 
{\bf Introduction} 
\vskip .1in

This is the first of a series of papers concerning   
what might be thought of as ``locally grouped spaces", in 
loose analogy with the locally ringed spaces of algebraic geometry. The spaces that we have in 
mind are simplicial sets that  generalize the simplicial sets that underly and 
determine  the classifying spaces of finite (or compact) groups.  If 
the analogy is pursued, then the role of ``structure sheaf" is provided by the ``fusion 
systems" associated with these spaces. Our approach here will be purely algebraic and combinatorial, 
so we will not be concerned with topological realizations.  
All of the groups to be considered will be finite; but a parallel series of papers representing  
some joint work with Alex Gonzalez, will considerably broaden the scope.

\vskip .1in 

Finite localities were introduced by the author in [Ch1], in order to give a positive solution to the 
question: Given a saturated fusion system $\ca F$ on a finite $p$-group, does there exist a 
``classifying space" for $\ca F$, and if so, is such a space unique up to isomorphism ?  
The solution that was given in [Ch1] was closely tied to the specific goal, and did not allow for 
a complete development of ideas. The aim here is to provide such a development.  In part, our aim is to 
supplement the theory of saturated fusion systems over a finite $p$-group. 
As part of that program we shall need to establish a sort of dictionary 
that will establish an equivalence between such notions as ``partial normal subgroup of a proper locality" 
and ``normal subsystem of a saturated fusion system". This is done in a seperate paper 
co-authored with Ellen Henke [ChHe].

The division into several papers closely parallels the extent to which fusion systems are drawn into the 
developing picture. This Part I can be characterized by its having no direct involvement with fusion systems, 
and by there being no mention in it of $p'$-elements or $p'$-subgroups of a group, other than in one 
application (see 4.12).

\vskip .1in 
Let $G$ be a group, and let $\bold W(G)$ be the free monoid on $G$. Thus, $\bold W(G)$ is the set of all 
words in the alphabet $G$, with the binary operation given by concatenation of words. The product 
$G\times G\to G$ extends, by generalized associativity, to a ``product" $\Pi:\bold W(G)\to G$, whereby a 
word $w=(g_1,\cdots,g_n)\in\bold W(G)$ is mapped to 
$g_1\cdots g_n$. The inversion map on $G$ induces an ``inversion" on $\bold W(G)$, sending  
$w$ to $(g_n\i,\cdots,g_1\i)$. In fact, one may easily replace the standard definition of ``group"
by a definition given in terms of $\Pi$ and the inversion on $\bold W(G)$. One obtains the notion of 
{\it partial group} by restricting the domain of $\Pi$ to a subset $\bold D$ of $\bold W(G)$, where 
$\bold D$, the product, and the inversion, are required to satisfy conditions (see definition 1.1) that 
preserve the outlines of the 
strictly group-theoretic setup. When one looks at things in this way, a group is simply a partial group 
$G$ having the property that $\bold D=\bold W(G)$. 

The notions of partial subgroup and homomorphism of partial groups immediately suggest themselves,  
and a partial subgroup of a partial group $\ca L$ may in fact be a group. We say that the partial group 
$\ca L$ is ``objective" (see definition 2.1) provided that the domain $\bold D$ of the product is 
determined in a certain way by a collection $\D$ of subgroups of $\ca L$ (the set of ``objects"), and 
provided that $\D$ has a certain ``closure" property. If also $\ca L$ is finite, and there exists 
$S\in\D$ such that $\D$ is a collection of subgroups of $S$, where $S$ is maximal in the set (partially 
ordered by inclusion) of $p$-subgroups of $\ca L$, then $(\ca L,\D,S)$ is a (finite) locality. 

The basic properties of partial groups, objective partial groups, and localities, will be derived in 
sections 1 and 2. We then begin in section 3 to consider partial normal subgroups of localities in  
detail.  One of the two key results in section 3 is the Frattini Lemma (3.11), which states 
that if $\ca N\norm\ca L$ is a partial normal subgroup, then $\ca L=N_{\ca L}(S\cap\ca N)\ca N$. 
The other is Stellmacher's splitting lemma (3.12), which leads to the partition of 
$\ca L$ into a 
collection of ``maximal cosets" of $\ca N$, and to a partial group structure on the set $\ca L/\ca N$ of 
maximal cosets. In section 4 it is shown that $\ca L/\ca N$ is in fact a locality, and we obtain 
versions of the first N\" other isomorphism theorem and of its familiar consequences. In particular, 
the notions ``partial normal subgroup" and ``kernel of a projection" turn out to be equivalent. This 
may be compared with the situation in the theory of saturated fusion systems, 
where it is known that no such equivalence exists.   

Section 5 concerns products of partial normal subgroups. The main result here (Theorem 5.1) has since 
been considerably strengthened by Ellen Henke [He], who shows that the product of any two partial 
normal subgroups of a locality is again a partial normal subgroup. The paper ends with a result 
(Proposition 5.5) which provides an application of essentially 
all of the concepts and results from all of the earlier sections, and which will play a role in Part III. 
\vskip .1in

Composition of mappings will most often be written from left to right, and mappings 
which are likely to be composed with others 
will be written to the right of their arguments. In particular, this entails that 
conjugation within a group $G$ be taken in the 
right-handed sense which is standard in finite group theory; so that $x^g=g\i xg$ for any $x,g\in G$. 

\vskip .1in 
The author extends his appreciation to Bernd Stellmacher for suggesting and proving the splitting lemma 
(3.12),  and for helpful suggestions regarding definition 1.1. Special thanks are due to Ellen Henke for a 
detailed list of corrections, and for her improvement in [He] on the results in section 5.

\vskip .2in 
\noindent 
{\bf Section 1: Partial groups} 
\vskip .1in 

The reader is asked to forget what a group is, and to trust that what was forgotten will soon be 
recovered. 

For any set $X$ write $\bold W(X)$ for the free monoid on $X$. Thus, an element 
of $\bold W(X)$ is a finite sequence of (or {\it word in}) the elements of $X$, and the multiplication in 
$\bold W(X)$ consists of concatenation of words, to be denoted $u\circ v$. The {\it length} of the word 
$(x_1,\cdots,x_n)$ is $n$. The {\it empty word} is the word $(\emptyset)$ of length 0. We make no 
distinction between $X$ and the set of words of length $1$. 

\definition {Definition 1.1} Let $\ca L$ be a non-empty set, let $\bold W=\bold W(\ca L)$ be the free monoid 
on $\ca L$, and let $\bold D$ be a subset of $\bold W$ such that: 
\roster 

\item "{(1)}" $\ca L\sub\bold D$ (i.e. $\bold D$ contains all words of length 1), and 
$$ 
u\circ v\in\bold D\implies u,v\in\bold D. 
$$
\endroster 
Notice that since $\ca L$ is non-empty, (1) implies that also the empty word is in $\bold D$. 

A mapping $\Pi:\bold D\to\ca L$ is a {\it product} if: 
\vskip .1in
\roster 

\item "{(2)}" $\Pi$ restricts to the identity map on $\ca L$, and 
\vskip .1in 

\item "{(3)}" $u\circ v\circ w\in\bold D\implies u\circ\Pi(v)\circ w\in\bold D$, and 
$\Pi(u\circ v\circ w)=\Pi(u\circ\Pi(v)\circ w)$. 
\vskip .1in 
\endroster 

An {\it inversion} on $\ca L$ consists of an involutory bijection $x\maps x\i$ on $\ca L$, 
together with the mapping $w\maps w\i$ on $\bold W$ given by 
$$ 
(x_1,\cdots,x_n)\maps(x_n\i,\cdots x_1\i). 
$$ 
We say that $\ca L$, with the product $\Pi:\bold D\to\ca L$ and inversion $(-)\i$, is a 
{\it partial group} if: 
\roster 

\item "{(4)}" $w\in\bold D\implies w\i\circ w\in\bold D$ and $\Pi(w\i\circ w)=\1$, 

\endroster 
where $\1$ denotes the image of the empty word under $\Pi$. Notice that (1) and (4) yield 
$w\i\in\bold D$ if $w\in\bold D$. As $(w\i)\i=w$, condition (4) is symmetric. 
\enddefinition 

\definition {Example 1.2} Let $\ca L$ be the $3$-element set $\{\1,a,b\}$ and let $\bold D$ be the subset of 
$\bold W(\ca L)$ consisting of all words $w$ such that the word obtained from $w$ by deleting all entries 
equal to $\1$ is an alternating string of $a$'s and $b$'s (of odd or even length and which, if non-empty, 
may begin either with $a$ or with $b$). Define $\Pi:\bold D\to\ca L$ by the formula: $\Pi(w)=\1$ if the number 
of $a$-entries in $w$ is equal to the number of $b$'s; $\Pi(w)=a$ if the number of $a$'s  
exceeds the number of $b$'s (necessarily by 1); and $\Pi(w)=b$ if the number of $b$'s 
exceeds the number of $a$'s. Define inversion on $\ca L$ by $\1\i=\1$, $a\i=b$, and $b\i=a$. 
It is then easy to check that $\ca L$, with these structures, is a partial group. In fact, $\ca L$ is 
the ``free partial group on one generator", as will be made clear in 1.12 below. 
\enddefinition 
 
It will be convenient to make the definition: a {\it group} is a partial group $\ca L$ in which 
$\bold W(\ca L)=\bold D$. In order to distinguish between this definition and the usual one we shall 
use the expression ``binary group" for a non-empty set $G$ with an associative binary operation, identity 
element, and inverses, in the usual sense. The following lemma shows that the distinction is subtle. 

\proclaim {Lemma 1.3} 
\roster  

\item "{(a)}" Let $G$ be a binary group and let $\Pi:\bold W(G)\to G$ be the ``multivariable 
 product" on $G$ given by $(g_1,\cdots,g_n)\maps g_1\cdots g_n$. Then $G$, together with $\Pi$ and 
the inversion in $G$, is a partial group, with $\bold D=\bold W(G)$. 

\item "{(b)}" Let $\ca L$ be a group; i.e. a partial group for which $\bold W(\ca L)=\bold D$. Then $\ca L$ 
is a binary group with respect to the operation given by restricting $\Pi$ to words of length 2, and 
with respect to the inversion in $\ca L$. Moreover, $\Pi$ is then the multivariable product on $\ca L$ 
defined as in (a). 

\endroster 
\endproclaim 

\demo {Proof} Point (a) is given by generalized associativity in the binary group $G$. Point 
(b) is a straightforward exercise, and is left to the reader. 
\qed 
\enddemo 

Here are a few basic consequences of definition 1.1. 

\proclaim {Lemma 1.4} Let $\ca L$ $($with $\bold D$, $\Pi$, and the inversion$)$ be a partial group. 
\roster 

\item "{(a)}" $\Pi$ is {\it $\bold D$-multiplicative}. That is, if $u\circ v$ is in $\bold D$ then 
the word $(\Pi(u),\Pi(v))$ of length 2 is in $\bold D$, and 
$$
\Pi(u\circ v)=\Pi(u)\Pi(v), 
$$
where ``$\Pi(u)\Pi(v)$" is an abbreviation for $\Pi((\Pi(u),\Pi(v))$.

\item "{(b)}" $\Pi$ is {\it $\bold D$-associative}. That is: 
$$
u\circ v\circ w\in\bold D\implies\Pi(u\circ v)\Pi(w)=\Pi(u)\Pi(v\circ w). 
$$ 

\item "{(c)}" If $u\circ v\in\bold D$ then $u\circ(\1)\circ v\in\bold D$ 
and $\Pi(u\circ(\1)\circ v)=\Pi(u\circ v)$. 

\item "{(d)}" If $u\circ v\in\bold D$ then both $u\i\circ u\circ v$ and 
$u\circ v\circ v\i$ are in $\bold D$, 
$\Pi(u\i\circ u\circ v)=\Pi(v)$, and $\Pi(u\circ v\circ v\i)=\Pi(u)$. 

\item "{(e)}" The cancellation rule: If both $u\circ v$ and $u\circ w$ are in $\bold D$, and 
$\Pi(u\circ v)=\Pi(u\circ w)$, then $\Pi(v)=\Pi(w)$ $($and similarly for 
right cancellation$)$. 

\item "{(f)}" If $u\in\bold D$ then $u\i\in\bold D$, and $\Pi(u\i)=\Pi(u)\i$. 
In particular, $\1\i=\1$. 

\item "{(g)}" The uncancellation rule: Suppose that both 
$u\circ v$ and $u\circ w$ are in $\bold D$ and that $\Pi(v)=\Pi(w)$. Then 
$\Pi(u\circ v)=\Pi(u\circ w)$. $($Similarly for right uncancellation.$)$ 

\endroster 
\endproclaim 

\demo {Proof} Let $u\circ v\in\bold D$. Then 1.1(3) applies to 
$(\nset)\circ u\circ v$ and yields 
$(\Pi(u))\circ v\in\bold D$ with $\Pi(u\circ v)=\Pi((\Pi(u))\circ v)$. Now 
apply 1.1(3) to $(\Pi(u))\circ v\circ(\nset)$, to obtain (a). 

Let $u\circ v\circ w\in\bold D$. Then $u\circ v$ and 
$w$ are in $\bold D$ by 1.1(1), and $\bold D$-multiplicativity yields 
$\Pi(u\circ v\circ w)=\Pi(u\circ v)\Pi(w)$. Similarly, 
$\Pi(u\circ v\circ w)=\Pi(u)\Pi(v\circ w)$, and (b) holds. 

Since $\1=\Pi(\nset)$, point (c) is immediate from 1.1(3). 

Let $u\circ v\in\bold D$. Then 
$v\i\circ u\i\circ u\circ v\in\bold D$ by 1.1(4), and then 
$u\i\circ u\circ v\in\bold D$ by 1.1(1). Multiplicativity then yields 
$$ 
\Pi(u\i\circ u\circ v)=\Pi(u\i\circ u)\Pi(v)=\1\Pi(v)= 
\Pi(\emptyset)\Pi(v)=\Pi(\emptyset\circ v)=\Pi(v). 
$$
As $(w\i)\i=w$ for any $w\in\bold W$, one obtains $w\circ w\i\in\bold D$ 
for any $w\in\bold D$, and $\Pi(w\circ w\i)=\1$. From this one easily completes 
the proof of (d). 

Now let $u\circ v$ and $u\circ w$ be in $\bold D$, with 
$\Pi(u\circ v)=\Pi(u\circ w)$. Then (d) (together with multiplicativity and 
associativity, which will not be explicitly mentioned hereafter) yield 
$$ 
\Pi(v)=\Pi(u\i\circ u\circ v)=\Pi(u\i)\Pi(u)\Pi(v)=\Pi(u\i)\Pi(u)\Pi(w)= 
\Pi(u\i\circ u\circ w)=\Pi(w), 
$$ 
and (e) holds. 

Let $u\in\bold D$. Then $u\circ u\i\in\bold D$, and then $\Pi(u)\Pi(u\i)=\1$. 
But also $(\Pi(u),\Pi(u)\i)\in\bold D$, and $\Pi(u)\Pi(u)\i=\1$. Now (f) 
follows by 1.1(2) and cancellation. 

Let $u,v,w$ be as in (g). Then $u\i\circ u\circ v$ and $u\i\circ u\circ w$ 
are in $\bold D$ by (d). By two applications of (d), 
$\Pi(u\i\circ u\circ v)=\Pi(v)=\Pi(w)=\Pi(u\i\circ u\circ w)$, so 
$\Pi(u\circ v)=\Pi(u\circ w)$ by (e), and (g) holds. 
\qed 
\enddemo

It will often be convenient to eliminate the symbol ``$\Pi$" and to speak of ``the product 
$g_1\cdots g_n$" 
instead of $\Pi(g_1,\cdots,g_n)$. More generally, if $\{X_i\}_{1\leq i\leq n}$ is a collection of subsets of 
$\ca L$ then the ``product set $X_1\cdots X_n$" is by definition the image under $\Pi$ of the set of words 
$(g_1,\cdots,g_n)\in\bold D$ such that $g_i\in X_i$ for all $i$. If $X_i=\{g_i\}$ is a singleton then we may 
write $g_i$ in place of $X_i$ in such a product. Thus, for example, the product 
$g\i Xg$ stands for the set of all $\Pi(g\i,x,g)$ with $(g\i,x,g)\in\bold D$, and with $x\in X$. 
\vskip .1in

{\it A Word of Urgent Warning}: In writing products in the above way one may be drawn into imagining that 
associativity holds in a stronger sense than that which is given by 1.4(b). This is an error that is to 
be avoided. For example one should not suppose, if $(f,g,h)\in\bold W$, and both $(f,g)$ and $(fg,h)$ are in 
$\bold D$, that $(f,g,h)$ is in $\bold D$. That is, it may be that ``the 
product $fgh$" is undefined, even though the product $(fg)h$ is defined. Of 
course, one is tempted to simply extend the domain $\bold D$ to include such 
triples $(f,g,h)$, and to ``define" the product $fgh$ to be $(fg)h$. The 
trouble is that it may also be the case that $gh$ and $f(gh)$ are defined, but that $(fg)h\neq f(gh)$.

\vskip .1in 
For $\ca L$ a partial group and $g\in\ca L$, write $\bold D(g)$ for the 
set of all $x\in\ca L$ such that the product $g\i xg$ is defined. There is then a mapping 
$$
c_g:\bold D(g)\to\ca L 
$$
given by $x\maps g\i xg$ (and called {\it conjugation by $g$}). Our 
preference is for right-hand notation for mappings, so we write 
$$ 
x\maps(x)c_g\quad\text{or}\quad x\maps x^g 
$$ 
for conjugation by $g$.

\vskip .1in 
The following result provides an illustration of the preceding notational conventions, 
and introduces a theme which will be developed further as we pass from partial groups 
to objective partial groups, localities, and (in Part III) regular localities.

\proclaim {Lemma 1.5} Let $\ca L$ be a partial group, and let $f,g\in\ca L$. 
\roster 

\item "{(a)}" Suppose that the products $fg$ and $gf$ are defined and that 
$fg=gf$. Suppose further that $f\in\bold D(g)$. Then $f^g=f$. 

\item "{(b)}" Suppose that $f\in\bold D(g)$, $g\in\bold D(f)$, and $f^g=f$. Then 
$fg=gf$ and $g^f=g$. 

\endroster 
\endproclaim 

\demo {Proof} (a): We are given $(f,g)\in\bold D$, so $(f\i,f,g)\in\bold D$ and $\Pi(f\i,f,g)=g$, 
by 1.4(d) and $\bold D$-associativity. We are given also $f\in\bold D(g)$ and $fg=gf$, so  
$$ 
f^g=\Pi(g\i,f,g)=\Pi((g\i,fg)=\Pi(g\i,gf)=\Pi(g\i,g,f)=f. 
$$
\vskip .1in 
\noindent 
(b): As $(g\i,f,g)\in\bold D$ we obtain $(f,g)\in\bold D$ from 1.1(1), and 
$(g,g\i,f,g)\in\bold D$ by 1.4(d). Then $\bold D$-associativity yields 
$fg=\Pi(g,g\i,f,g)=gf^g$. As $f^g=f$ by hypothesis, we obtain $fg=gf$. Finally, since 
$(f\i,f,g)$ and $(f\i,g,f)$ are in $\bold D$ the uncancellation rule yields 
$f\i fg=f\i gf$, and so $g^f=g$. 
\qed 
\enddemo

\definition {Notation} From now on, in any given partial group $\ca L$, usage of the symbol ``$x^g$" shall 
be taken to imply $x\in\bold D(g)$. More generally, for $X$ a subset of $\ca L$ and $g\in\ca L$, usage of 
``$X^g$" shall be taken to mean that $X\sub\bold D(g)$; whereupon $X^g$ is by definition the set of all 
$x^g$ with $x\in X$. 
\enddefinition

At this early point, and in the context of arbitrary partial groups, one can say very little about the 
maps $c_g$. The cancellation rule 1.4(e) implies 
that each $c_g$ is injective, but beyond that the following lemma may be the best that can be obtained.

\proclaim {Lemma 1.6} Let $\ca L$ be a partial group and let $g\in\ca L$. Then the following hold. 
\roster

\item "{(a)}" $\1\in\bold D(g)$ and ${\1}^g=\1$. 

\item "{(b)}" $\bold D(g)$ is closed under inversion, and $(x\i)^g=(x^g)\i$ for all $x\in\bold D(g)$. 

\item "{(c)}" $c_g$ is a bijection $\bold D(g)\to\bold D(g\i)$, and $c_{g\i}=(c_g)\i$. 

\item "{(d)}" $\ca L=\bold D(\1)$, and $x^{\1}=x$ for each $x\in\ca L$. 

\endroster 
\endproclaim 

\demo {Proof} By 1.1(4), $g\circ\nset\circ g\i=g\circ g\i\in\bold D$, so $\1\in\bold D(g)$ and then 
$\1^g=\1$ by 1.4(c). Thus (a) holds. Now let $x\in\bold D(g)$ and set $w=(g\i,x,g)$. Then $w\in\bold D$, 
and $w\i=(g\i,x\i,g)$ by definition in 1.1. Then 1.1(4) yields $w\i\circ w\in\bold D$, and so 
$w\i\in\bold D$ by 1.1(1). This shows that $\bold D(g)$ is closed under inversion. Also, 1.1(4) yields 
$\1=\Pi(w\i\circ w)=(x\i)^g x^g$, and then $(x\i)^g=(x^g)^{\i}$ by 1.4(f). This completes the proof of (b). 

As $w\in\bold D$, 1.4(d) implies that $g\circ w$ and then $g\circ w\circ g\i$ are in $\bold D$. Now 1.1(3) 
and two applications of 1.4(d) yield 
$$ 
gx^gg\i=\Pi(g,g\i,x,g,g\i)=\Pi((g,g\i,x)\circ g\circ g\i)=\Pi(g,g\i,x)=x. 
$$ 
Thus $x^g\in\bold D(g\i)$ with $(x^g)^{g\i}=x$, and thus (c) holds. 

Finally, $\1=\1\i$ by 1.4(f), and $\nset\circ x\circ\nset=x\in\bold D$ for any $x\in\ca L$, proving (d). 
\qed 
\enddemo 

\definition {Definition 1.7} Let $\ca L$ be a partial group and let $\ca H$ be a non-empty subset of 
$\ca L$. Then $\ca H$ is a {\it partial subgroup} of $\ca L$ (denoted $\ca H\leq\ca L$) if $\ca H$ is 
closed under inversion ($g\in\ca H$ implies $g\i\in\ca H$) and closed with respect to products. The latter 
condition means, of course, that $\Pi(w)\in\ca H$ whenever $w\in\bold W(\ca H)\cap\bold D$. A partial subgroup 
$\ca N$ of $\ca L$ is a partial {\it normal} subgroup of $\ca L$ (denoted $\ca N\norm\ca L$) if $x^g\in\ca N$ 
for all $x\in\ca N$ and all $g\in\ca L$ for which $x\in\bold D(g)$. 
We say that $\ca H$ is a {\it subgroup} of $\ca L$ if $\ca H\leq\ca L$ and $\bold W(\ca H)\sub\bold D$. 

An equivalent way to state the condition for normality, which relies on the notational convention introduced 
above for interpreting product sets $XYZ$, is to say that the 
partial subgroup $\ca N$ of $\ca L$ is normal in $\ca L$ if $g\i\ca Ng\sub\ca N$ for all $g\in\ca L$. 
\enddefinition

We leave it to the reader to check that if $\ca H\leq\ca L$ then $\ca H$ is indeed a partial group, with 
$\bold D(\ca H)=\bold W(\ca H)\cap\bold D(\ca L)$.

\proclaim {Lemma 1.8} Let $\ca H$ and $\ca K$ be partial subgroups of a partial group $\ca L$, and let 
$\{\ca H_i\}_{i\in I}$ be a set of partial subgroups of $\ca L$. 
\roster 

\item "{(a)}" Each partial subgroup of $\ca H$ is a partial subgroup of $\ca L$. 

\item "{(b)}" Each partial subgroup of $\ca L$ which is contained in $\ca H$ is a partial subgroup of $\ca H$. 

\item "{(c)}" If $\ca H$ is a subgroup of $\ca L$ then $\ca H\cap\ca K$ is a subgroup of both $\ca H$ and 
$\ca K$.  

\item "{(d)}" Suppose $\ca K\norm\ca L$. Then $\ca H\cap\ca K\norm\ca H$. 
Moreover, $\ca H\cap\ca K$ is a normal subgroup of $\ca H$ if $\ca H$ is a subgroup of $\ca L$. 

\item "{(e)}" $\bigcap\{\ca H_i\mid i\in I\}$ is a partial subgroup of $\ca L$, and is 
a partial normal subgroup of $\ca L$ if $\ca H_i\norm\ca L$ for all $i$.  

\endroster  
\endproclaim 

\demo {Proof} One observes that in all of the points (a) through (e) the requisite closure with respect to 
inversion obtains. Thus, we need only be concerned with products. 
\vskip .1in 
\noindent 
(a) Let $\ca E\leq\ca H$ be a partial subgroup of $\ca H$. Then 
$$ 
\bold D(\ca E)=\bold W(E)\cap\bold D(\ca H)=\bold W(E)\cap(\bold W(\ca H)\cap\bold D(\ca L)=
\bold W(E)\cap\bold D(\ca L), 
$$ 
and (a) follows. 
\vskip .1in 
\noindent 
(b) Suppose $\ca K\sub\ca H$ and let $w\in\bold W(\ca K)\cap\bold D(\ca H)$. As 
$\bold D(\ca H)\leq\bold D(\ca L)$, and since $\ca K\leq\ca L$ by hypothesis, we obtain $\Pi(w)\in\ca K$. 
\vskip .1in 
\noindent 
(c) Assuming now that $\ca H$ is a subgroup of $\ca L$, we have $\bold W(\ca H)\sub\bold D(\ca L)$, and 
then $\bold D(\ca H\cap\ca K)\sub\bold D(\ca H)\cap\bold D(\ca K)$, so that $\ca H\cap\ca K$ is a subgroup of 
both $\ca H$ and $\ca K$. 
\vskip .1in 
\noindent 
(d) Let $\ca K\norm\ca L$ and let $x\in\ca H\cap\ca K$ and $h\in\ca H$ with $(h\i,x,h)\in\bold D(\ca H)$. 
Then $(h\i,x,h)\in\bold D(\ca L)$, and $x^h\in\ca K$. As $\ca H\leq\ca L$ we have also $x^h\in\ca H$, and 
so $\ca H\cap\ca K\norm\ca H$. Now suppose further that $\ca H$ is a subgroup of $\ca L$. That is, assume 
that $\bold W(\ca H)\sub\bold D(\ca L)$. Then $\bold W(\ca H\cap\ca K)\sub\bold D(\ca L)$, hence 
$\ca H\cap\ca K$ is a subgroup of $\ca H$, and evidently a normal subgroup.   
\vskip .1in 
\noindent 
(e) Set $\ca X=\bigcap\{\ca H_i\}_{i\in I}$. Then $\Pi(w)\in\ca X$ for all 
$w\in\bold W(\ca X)\cap\bold D(\ca L)$, 
and so $\ca X\leq\ca L$. The last part of (e) may be left to the reader. 
\qed 
\enddemo

For any subset $X$ of a partial group $\ca L$ define the partial subgroup $\<X\>$ of $\ca L$ 
{\it generated by $X$} to be the intersection of the set of all partial subgroups of $\ca L$ containing $X$. 
Then $\<X\>$ is itself a partial subgroup of $\ca L$ by 1.8(e).

\proclaim {Lemma 1.9} Let $X$ be a subset of $\ca L$ such that $X$ is closed under inversion. Set 
$X_0=X$ and recursively define $X_n$ for $n>0$ by 
$$ 
X_n=\{\Pi(w)\mid w\in\bold W(X_{n-1})\cap\bold D\}. 
$$ 
Then $\<X\>=\bigcup\{X_n\}_{n\geq 0}$. 
\endproclaim 

\demo {Proof} Let $Y$ be the union of the sets $X_i$. Each $X_i$ is closed under inversion by 1.4(f), and 
$Y\neq\nset$ since $\1=\Pi(\nset)$. Since $Y$ is closed under products, by construction, we get $Y\leq\<X\>$, 
and then $Y=\<X\>$ by the definition of $\<X\>$. 
\qed 
\enddemo

\proclaim {Lemma 1.10 (Dedekind Lemma)} Let $\ca H$, $\ca K$, and $\ca A$ be partial subgroups of a 
partial group $\ca L$, and assume $\ca A$ is a subset of $\ca H\ca K$.  
\roster 

\item "{(a)}" If $\ca K\leq\ca A$ then $\ca A=(\ca A\cap\ca H)\ca K$. 

\item "{(b)}" If $\ca H\leq\ca A$ then $\ca A=\ca H(\ca A\cap\ca K)$.

\endroster  
\endproclaim 

\demo {Proof} The proof is identical to the proof for binary groups, and is left to the reader. 
\qed 
\enddemo 

\definition {Definition 1.11} Let $\ca L$ and $\ca L'$ be partial groups, let $\b:\ca L\to\ca L'$ be a 
mapping, and let $\b^*:\bold W\to\bold W'$ be the induced mapping of free monoids. Then $\b$ is a 
{\it homomorphism $($of partial groups$)$} if: 
\roster 

\item "{(H1)}" $\bold D\b^*\sub\bold D'$, and 

\item "{(H2)}" $(\Pi(w))\b=\Pi'(w\b^*)$ for all $w\in\bold D$. 

\endroster 
The {\it kernel} of $\b$ is the set $Ker(\b)$ of all $g\in\ca L$ such that $g\b=\1'$. We say that $\b$ is an 
{\it isomorphism} if there exists a homomorphism $\b':\ca L'\to\ca L$ such that $\b\circ\b'$ and 
$\b'\circ\b$ are identity mappings. (Equivalently, $\b$ is an isomorphism if $\b$ is bijective and 
$\bold D\b=\bold D'$.) 
\enddefinition 

\definition {Example 1.12} Let $\ca L=\{\1,a,b\}$ be the partial group from example 1.2, let $\ca L'$ be 
any partial group, and let $x\in\ca L'$. Then the mapping $\b:\ca L\to\ca L'$ given by 
$$ 
\1\maps\1', \quad a\maps x, \quad  b\maps x\i. 
$$ 
is a homomorphism. In fact, $\b$ is the unique homomorphism $\ca L\to\ca L'$ which maps $a$ to $x$, 
by the following lemma. Thus, $\ca L$ is the (unique up to a unique invertible homomorphism) 
{\it free partial group} on one generator. Free partial groups in general can be obtained as ``free 
products" of copies of $\ca L$ (see Appendix A). 
\enddefinition

\proclaim {Lemma 1.13} Let $\b:\ca L\to\ca L'$ be a homomorphism of partial 
groups. Then $\1\b=\1'$, and $(g\i)\b=(g\b)\i$ for all $g\in\ca L$. 
\endproclaim 

\demo {Proof} Since $\1 \1=\1$, (H1) and (H2) yield $\1\b=(\1\1)\b=(\1\b)(\1\b)$, and then $\1\b=\1'$ 
by left or right cancellation. Since $(g,g\i)\in\bold D$ for any $g\in\ca L$ by 1.4(d),
(H1) yields $(g\b,(g\i)\b)\in\bold D'$, and then
$\1\b=(gg\i)\b=(g\b)((g\i)\b)$ by (H2). As $\1\b=\1'=(g\b)(g\b)\i$,
left cancellation yields $(g\i)\b=(g\b)\i$.
\qed
\enddemo

\proclaim {Lemma 1.14} Let $\b:\ca L\to\ca L'$ be a homomorphism of partial groups, and set $\ca N=Ker(\b)$. 
Then $\ca N$ is a partial normal subgroup of $\ca L$.  
\endproclaim 

\demo {Proof} By 1.13 $\ca N$ is closed under inversion. For $w$ in $\bold W(\ca N)\cap\bold D$ the map 
$\b^*:\bold W\to\bold W'$ sends $w$ to a word of the form $(\1',\cdots,\1')$. Then $\Pi'(w\b^*)=\1'$, 
and thus $\Pi(w)\in\ca N$ and $\ca N$ is a partial subgroup of $\ca L$. Now let 
$f\in\ca L$ and let $g\in\ca N\cap\bold D(f)$. Then 
$$
(f\i,g,f)\b^*=((f\b)\i,{\1}',f\b) \quad\text{(by 1.13)}, 
$$
so that 
$$
(g^f)\b=\Pi'((f\i,g,f)\b^*)=\Pi'(f\b)\i,\1',f\b)=\1'. 
$$ 
Thus $\ca N\norm\ca L$. 
\qed 
\enddemo

It will be shown later (cf. 4.6) that partial normal subgroups of ``localities" are always 
kernels of homomorphisms.

\proclaim {Lemma 1.15} Let $\b:\ca L\to\ca L'$ be a homomorphism of partial groups and let $M$ be a 
subgroup of $\ca L$. Then $M\b$ is a subgroup of $\ca L'$. 
\qed 
\endproclaim 

\demo {Proof} We are given $\bold W(M)\sub\bold D(\ca L)$, so $\b^*$ maps $\bold W(M)$ into 
$\bold D(\ca L')$. (Note, however, example 1.12.) 
\qed 
\enddemo

\proclaim {Lemma 1.16} Let $G$ and $G'$ be groups (and hence also binary groups in the sense of 1.3). 
A map $\a:G\to G'$ is a homomorphism of partial groups if and only if $\a$ is a homomorphism of 
binary groups. 
\endproclaim  

\demo {Proof} We leave to the reader the proof that if $\a$ is a homomorphism of partial groups then 
$\a$ is a homomorphism of binary groups. Now suppose that $\a$ is a homomorphism of binary groups. 
As $\bold W(G)=\bold D(G)$ (and similarly for $G'$, it is immediate that $\a^*$ maps $\bold D(G)$ into 
$\bold D(G)$. Assume that $\a$ is not a 
homomorphism of partial groups and let $w\in\bold D(G)$ be of minimal length subject to 
$\Pi'(w\a^*)\neq(\Pi(w))\a$. Then $n>1$ and we can write $w=u\circ v$ with both $u$ and $v$ non-empty. Then 
$$ 
\Pi'(w\a^*)=\Pi'(u\a^*\circ v\a^*)=\Pi'(u\a^*)\Pi'(v\a^*)=((\Pi(u))\a)((\Pi(v))\a)=(\Pi(u)\Pi(v))\a, 
$$ 
as $\a$ is a homomorphism of binary groups. Since $(\Pi(u)\Pi(v))\a=(\Pi(w))\a$, the proof is complete. 
\qed 
\enddemo

\vskip .2in 
\noindent 
{\bf Section 2: Objective partial groups and localities} 
\vskip .1in 

Recall the convention: if $X$ is a subset of the partial group $\ca L$, and $g\in\ca L$, then any statement 
involving the expression ``$X^g$" is to be understood as carrying the assumption that $X\sub\bold D(g)$. 
Thus, the statement ``$X^g=Y$" means: $(g\i,x,g)\in\bold D$ for all $x\in X$, and $Y$ is the set of 
products $g\i xg$ with $x\in X$. 

\definition {Definition 2.1} Let $\ca L$ be a partial group. For any collection $\D$ of subgroups of 
$\ca L$ define $\bold D_\D$ to be the set of all $w=(g_1,\cdots,g_n)\in\bold W(\ca L)$ such that:
\roster

\item "{(*)}" there exists $(X_0,\cdots,X_n)\in\bold W(\D)$ with
$(X_{i-1})^{g_i}=X_i$ for all $i$ ($1\leq i\leq n$).

\endroster 
Then $\ca L$ is {\it objective} if there exists a set $\D$ of subgroups of $\ca L$ 
such that the following two conditions hold.
\roster

\item "{(O1)}" $\bold D=\bold D_{\D}$.

\item "{(O2)}" Whenever $X$ and $Y$ are in $\D$ and $g\in\ca L$ such that $X^g$ is a 
subgroup of $Y$, then every subgroup of $Y$ containing $X^g$ is in $\D$. 

\endroster 
We say also that $\D$ is a set of {\it objects}) for $\ca L$, if (O1) and (O2) hold.  
\enddefinition

It will often be convenient to somewhat over-emphasize the role of $\D$ in the above definition by saying 
that ``$(\ca L,\D)$ is an objective partial group". What is meant by this is that $\ca L$ is an objective 
partial group and that $\D$ is a set (there will often be more than one) of objects for $\ca L$. 

We mention that the condition (O2) requires more than that $X^g$ be a sub{\it set} of $Y$, in order to 
conclude that overgroups of $X$ in $Y$ are objects. This is a non-vacuous distinction, since the 
conjugation map $c_g:X\to X^g$ need not send $X$ to a subgroup 
of $\ca L$, in a general partial group. 

\definition {Example 2.2} Let $G$ be a group, let $S$ be a subgroup of $G$, and let $\D$ be 
a collection of subgroups of $S$ such that $S\in\D$. Assume that $\D$ satisfies (O2). That is, assume that 
$Y\in\D$ for every 
subgroup $Y$ of $S$ such that $X^g\leq Y$ for some $X\in\D$ and some $g\in G$. Let $\ca L$ 
be the set of all $g\in G$ such that $S\cap S^g\in\D$, and let $\bold D$ be the subset 
$\bold D_\D$ of $\bold W(\ca L)$. Then $\ca L$ is a partial group (via the multivariable 
product in $G$ and the inversion in $G$), and $(\ca L,\D)$ is an objective partial group. Specifically: 
\roster 

\item "{(a)}" If $\D=\{S\}$ then $\ca L=N_G(S)$ (and so $\ca L$ is a group in this case). 

\item "{(b)}" Take $G=O_4^+(2)$. Thus, $G$ is a semidirect product $V\rtimes S$ where $V$ is elementary 
abelian of order 9 and $S$ is a dihedral group of order 8 acting faithfully on $V$. Let $\D$ be the set of 
all non-identity subgroups of $S$. One may check that $S\cap S^g\in\D$ for all $g\in G$, and hence 
$\ca L=G$ (as {\it sets}). But $\ca L$ is not a group, as $\bold D_\D\neq\bold W(G)$. 

\item "{(c)}" Take $G=GL_3(2)$, $S\in Syl_2(G)$, and let $M_1$ and $M_2$ be the two maximal subgroups of $G$ 
containing $S$. Set $P_i=O_2(M_i)$ and set $\D=\{S,P_1,P_2\}$. Then $\ca L=M_1\cup M_2$ (in fact the ``free 
amalgamated product" of $M_1$ and $M_2$ over $S$ in the category of partial groups). On 
the other hand, if $\D$ is taken to be the set of all non-identity subgroups of $S$ then $\ca L$ is 
somewhat more complicated. Its underlying set is $M_1M_2\cup M_2M_1$. 

\endroster 
\enddefinition

\vskip .1in 
In an objective partial group $(\ca L,\D)$ we say that the word $w=(g_1,\cdots,g_n)$ is 
{\it in $\bold D$ via $(X_0,\cdots,X_n)$} if the condition (*) in 2.1 applies specifically 
to $w$ and $(X_0,\cdots,X_n)$. We may also say, more simply, that {\it $w$ is in $\bold D$ 
via $X_0$}, since the sequence $(X_0,\cdots,X_n)$ is determined by $w$ and $X_0$. 

\vskip .1in 
For any partial group $\ca L$ and subgroups $X,Y$ of $\ca L$, set   
$$ 
N_{\ca L}(X,Y)=\{g\in\ca L\mid X\sub\bold D(g),\ X^g\sub Y\},   
$$ 
and set 
$$ 
N_{\ca L}(X)=\{g\in\ca L\mid X^g=X\}. 
$$ 

\proclaim {Lemma 2.3} Let $(\ca L,\D)$ be an objective partial group. 
\roster 

\item "{(a)}" $N_{\ca L}(X)$ is a subgroup of $\ca L$ for each $X\in\D$. 

\item "{(b)}" Let $g\in\ca L$ and let $X\in\D$ with $Y:=X^g\in\D$. Then $N_{\ca L}(X)\sub\bold D(g)$, and 
$$
c_g:N_{\ca L}(X)\to N_{\ca L}(Y) 
$$
is an isomorphism of groups. More generally:

\item "{(c)}" Let $w=(g_1,\cdots,g_n)\in\bold D$ via $(X_0,\cdots,X_n)$. 
Then 
$$
c_{g_1}\circ\cdots\circ c_{g_n}=c_{\Pi(w)} 
$$
as isomorphisms from $N_G(X_0)$ to $N_G(X_n)$. 

\endroster 
\endproclaim 

\demo {Proof}(a) Let $X\in\D$ and let $u\in\bold W(N_{\ca L}(X))$. Then $u\in\bold D$ 
via $X$, $\1\in N_{\ca L}(X)$ (1.6(d)), and $N_{\ca L}(X)\i=N_{\ca L}(X)$ (1.6(c)). 
\vskip .1in 
\noindent 
(b) Let $x,y\in N_{\ca L}(X)$ and set $v=(g\i,x,g,g\i,y,g)$. Then $v\in\bold D$ via $Y$, 
and then $\Pi(v)=(xy)^g=x^gy^g$ (using points (a) and (b) of 1.4). Thus, the 
conjugation map $c_g:N_{\ca L}(X)\to N_{\ca L}(Y)$ is a homomorphism of binary groups  
(see 1.3), and hence a homomorphism of partial groups (1.16). Since $c_{g\i}=c_g\i$ by 1.6(c), 
$c_g$ is an isomorphism of groups. 
\vskip .1in 
\noindent 
(c) Let $x\in N_{\ca L}(X_0)$, set $u_x=w\i\circ(x)\circ w$, and observe that $u_x\in\bold D$ 
via $X_n$. Then $\Pi(u_x)$ can be written as $(\cdots(x)^{g_1}\cdots)^{g_n}$, and this yields (c). 
\qed 
\enddemo

The next lemma provides two basic computational tools.

\proclaim {Lemma 2.4} Let $(\ca L,\D)$ be an objective partial group. 
\roster 

\item "{(a)}"  Let $(a,b,c)\in\bold D$ and set $d=abc$. Then $bc=a\i d$ and $ab=dc\i$ 
(and all of these products are defined). 

\item "{(b)}" Let $(f,g)\in\bold D$ and let $X\in\D$. Suppose that both 
$X^f$ and $X^{fg}$ are in $\D$. Then $X^{fg}=(X^f)^g$. 

\endroster 
\endproclaim 

\demo {Proof} Point (a) is a fact concerning partial groups in general, and is immediate from 1.4(c). 
Now consider the setup in (b). As $(f,g)\in\D$ we have also $(f\i,f,g)\in\D$, and $g=\Pi(f\i,f,g)=f\i(fg)$. 
Now observe that $(f\i,fg)\in\bold D$ via $P^f$, and apply 2.3(c) to obtain 
$P^{fg}=((P^f)^{f\i})^{fg}=(P^f)^g$.  
\qed 
\enddemo 

The following result is a version of lemma 1.5(b) for objective partial groups. 
The hypothesis is weaker than that of 1.5(b), and the conclusion is stronger.

\proclaim {Lemma 2.5} Let $(\ca L,\D)$ be an objective partial group and let $f,g\in\ca L$. 
Suppose that $f^g=f$. Then $fg=gf$ and $g^f=g$. 
\endproclaim 

\demo {Proof} Suppose that $(g\i,f,g)\in\bold D$ via $(P_0,P_1,P_2,P_3)$. One then has the following 
commutative square of conjugation maps, in which the arrows are labeled by elements that perform the 
conjugation.
$$ 
\CD 
P_2@>g>>P_3 \\ 
@AfAA  @AA{f^g}A \\ 
P_1@>>g>P_0 
\endCD 
$$ 
Now assume that $f^g=f$. Since any of the arrows in the diagram may be reversed using 1.6(c), one reads 
off that $(f\i,g,f)\in\bold D$ via $P_2$. Then $g^f=g$ and $fg=gf$ by 1.5(b). 
\qed 
\enddemo 

The following result and its corollary are fundamental to the entire enterprise. The proof given here is 
due to Bernd Stellmacher. 

\proclaim {Proposition 2.6} Let $(\ca L,\D)$ be an objective partial group. Suppose that $\D$ is a 
collection of subgroups of some $S\in\D$. For each $g\in\ca L$ define 
$S_g$ to be the set of all $x\in\bold D(g)\cap S$ such that $x^g\in S$. Then: 
\roster 

\item "{(a)}" $S_g\in\D$. In particular, $S_g$ is a subgroup of $S$. 

\item "{(b)}" The conjugation map $c_g:S_g\to(S_g)^g$ is an isomorphism of groups, and 
$S_{g\i}=(S_g)^g$. 

\item "{(c)}" $P^g$ is defined and is a subgroup of $S$, for every subgroup $P$ of $S_g$. In 
particular, $P^g\in\D$ for any $P\in\D$ with $P\leq S_g$. 

\endroster 
\endproclaim 

\demo {Proof} Fix $g\in\ca L$. Then the word $(g)$ of length 1 is in $\bold D$ by 1.1(2), and since  
$\bold D=\bold D_\D$ by (O1) there exists $X\in\D$ such that $Y:=X^g\in\D$. Let $a\in S_g$ 
and set $b=a^g$. Then $X^a$ and $X^b$ are subgroups of $S$ (as $a,b\in S$), so $X^a$ and 
$Y^b$ are in $\D$ by (O2). Then $(a\i,g,b)\in\bold D$ via $X^a$, so also $(g,b)\in\bold D$. Also  
$(a,g)\in\bold D$ via $X^{a\i}$. Since $g\i ag=b$ we get $ag=gb$ by cancellation, and hence 
$$ 
a\i gb=a\i(gb)=a\i(ag)=(a\i a)g=g 
$$ 
by $\bold D$-associativity. Since $a\i gb$ conjugates $X^a$ to $Y^b$, we 
draw the following conclusion. 
\roster 

\item "{(1)}" $X^a\leq S_g$ and $(X^a)^g\in\D$ for all $a\in S_g$. 

\endroster 
Now let $c,d\in S_g$. Then (1) shows that both $X^c$ and $X^{cd}$ are members  
of $\D$ which are conjugated to members of $\D$ by $g$. Setting 
$w=(g\i,c,g,g\i,d,g)$, we conclude (by following $X^g$ along the chain of 
conjugations given by $w$) that $w\in\bold D$ via $X^g$. Then $\bold D$-associativity  
yields 
$$ 
\Pi(w)=(cd)^g=c^g d^g.  \tag2
$$ 
Since $c^g$ and $d^g$ are in $S$, we conclude that $cd\in S_g$. Since $S_g$ is closed under inversion 
by 1.6(b), $S_g$ is a subgroup of $S$. As $X\leq S_g\leq S$, where $X$ and $S$ are in $\D$, (02) now 
yields $S_g\in\D$. Thus (a) holds. 

Since $c_{g\i}=(c_g)\i$ by 1.6(c), it follows that $S_{g\i}=(S_g)^g$. Points (b) and (c) are then 
immediate from (a) and 2.3(b). 
\qed 
\enddemo

\proclaim {Corollary 2.7} Assume the hypothesis of 2.6, let $w=(g_1,\cdots,g_n)\in\bold W(\ca L)$, and 
define $S_w$ to be the set of all $x\in S$ such that, for all $k$ with $1\leq k\leq n$, the composition 
$c_{g_1}\circ\cdots\circ c_{g_k}$ is defined on $x$ and maps $x$ into $S$. 
Then $S_w$ is a subgroup of $S$, and $S_w\in\D$ if and only if $w\in\bold D$. 
\endproclaim 

\demo {Proof} Let $x_0,y_0\in S_w$ and let $\s=(x_0,\cdots,x_n)$ and $\tau=(y_0,\cdots,y_n)$ be the 
corresponding sequences of elements of $S$, obtained from $x_0$ and $y_0$ via the sequence of maps 
$c_{g_1}\circ\cdots\circ c_{g_k}$. Set $S_j=S_{g_j}$ $(1\leq j\leq n)$. 
Then $x_{i-1}$ and $y_{i-1}$ are elements of $S_{i-1}$, and so $x_{i-1}y_{i-1}\in S_{i-1}$ by 2.6(a). As 
$c_{g_i}$ restricts to a homomorphism on $S_{i-1}$ (see 2.3(b)) it follows that $x_iy_i\in S_i$. 
Thus $S_w$ is closed under the binary product in $S$. That $S_w$ is closed under inversion is given 
by 1.6(c), so $S_w$ is a subgroup of $S$. 

Suppose that $S_w\in\D$, set $P_0=S_w$, and recursively define $P_k$ for $0\<k\leq n$ by 
$P_k=(P_{k-1})^{g_k}$. Then $P_k$ is a subgroup of $S$ by 2.6(c) and induction on $k$. Then 
(O2) yields $P_k\in\D$, and so $w\in\bold D$ by (O1). Conversely, if $w\in\bold D$ then 
(O1) shows that $P\leq S_w$ for some $P\in\D$, and then $S_w\in\D$ by (O2). 
\qed 
\enddemo

Henceforth our focus will be on {\it finite} objective partial groups 
of a certain kind. Here is the main definition.

\definition {Definition 2.8} Let $p$ be a prime and let $\ca L$ be a finite partial group. Then $\ca L$ is 
a {\it locality} if there exists a $p$-subgroup $S$ of $\ca L$ and a set $\D$ of subgroups of $S$ such 
that: 
\roster 

\item "{(L1)}" $(\ca L,\D)$ is objective.  

\item "{(L2)}" $S$ is in $\D$, and $S$ is maximal in the set (partially ordered by inclusion) of 
$p$-subgroups of $\ca L$. 

\endroster 
\enddefinition 

As with definition 2.1, we shall tend to over-emphasize the roles of $S$ and of $\D$ by saying that 
``$(\ca L,\D,S)$ is a locality" when, strictly speaking, we mean only that $\ca L$ is a locality and that 
$S$ and $\D$ fulfill the conditions (L1) and (L2). (The extent to which $S$ and $\D$ are 
determined by $\ca L$ is explored in 2.16 and 2.17 below.)  

Notice that if $\ca L$ is a locality then the hypothesis of 2.6 and of 2.7 is fulfilled, and 
we may therefore speak of the subgroups $S_g$ and $S_w$ of $S$ for any $g\in\ca L$ and 
any $w\in\bold W(\ca L)$.

\proclaim {Lemma 2.9} Let $(\ca L,\D,S)$ be a locality, let $(x_1,\cdots,x_n)\in\bold D$, and let 
$(g_0,\cdots,g_n)\in\bold W(N_{\ca L}(S))$. Then 
$$ 
(g_0,x_1,g_1\i,g_1,x_2,\cdots,x_{n_1},g_{n-1}\i,g_{n-1},x_n,g_n)\in\bold D. \tag*
$$ 
In particular, conjugation by $g\in N_{\ca L}(S)$ is an automorphism of the partial group $\ca L$. 
\endproclaim 

\demo {Proof} Set $w=(x_1,\cdots,x_n)$ and set $P=(S_w)^{g_0\i}$. Then the word displayed in (*) is in 
$\bold D$ via $P$. 
\qed 
\enddemo

\proclaim {Lemma 2.10} Let $(\ca L,\D,S)$ be a locality and let $P\in\D$. Then there exists $g\in\ca L$ 
such that $N_S(P)\leq S_g$ and such that $N_S(P^g)$ is a Sylow $p$-subgroup of $N_{\ca L}(P^g)$. 
\endproclaim  

\demo {Proof} Observe first of all that, by (L2) the lemma holds for $P=S$ and $g=\1$.   
Among all $P$ for which the lemma fails 
to hold, choose $P$ so that first $|P|$ and then $|N_S(P)|$ is as large as possible. Set $R=N_S(P)$ 
and let $R^*$ be a Sylow $p$-subgroup of $N_{\ca L}(P)$ containing $R$. 
Then $R<R^*$ (proper subgroup), and then also $R<N_{R^*}(R)$. We have $P\neq R$ as $P\neq S$,  
and then the maximality of $|P|$ yields the existence of 
an element $f\in N_{\ca L}(R,S)$ with $N_S(R^f)\in Syl_p(N_{\ca L}(R^f))$. 

By 2.3(b) $f$-conjugation induces an isomorphism 
$$ 
N_{\ca L}(R)@>c_f>>N_{\ca L}(R^f).  
$$ 
By Sylow's theorem there exists $x\in N_{\ca L}(R^f)$ such that $(N_{R^*}(R)^f)^x\leq N_S(R^f)$. Here 
$(f,x)\in\bold D$ via $R$, so 2.3(c) yields $(N_{R^*}(R)^f)^x=N_{R^*}(R)^{fx}$. Thus, 
by replacing $f$ with $fx$, we may assume that $f$ was chosen to begin with so that 
$N_{R^*}(R)^f\leq N_S(R^f)$. Since $R^*$ normalizes $P$ and $c_f$ is an 
isomorphism, it follows that $N_{R^*}(R)^f$ normalizes $P^f$, and thus 
$|N_S(P^f)|>|N_S(P)|$. The maximality of $|N_S(P)|$ in the choice of $P$ then implies 
that $P^f$ is not a counter-example to the lemma. Set $Q=P^f$. Thus there exists 
$h\in N_{\ca L}(N_S(Q),S)$ such that $N_S(Q^h)$ is a Sylow subgroup of $N_{\ca L}(Q^h)$. 
Here $(f,h)\in\bold D$ via $R$, so $Q^h=P^g$ where $g=fh$, so $P$ is not a counter-example to the 
lemma, and therefore no counter-example exists.  
\qed 
\enddemo

\proclaim {Proposition 2.11} Let $(\ca L,\D,S)$ be a locality and let $H$ be a subgroup of $\ca L$. 
\roster 

\item "{(a)}" There exists an object $P\in\D$ such that $H\leq N_{\ca L}(P)$. Indeed, there exists 
a unique largest such $P$. $($All subgroups are ``local subgroups".$)$. 

\item "{(b)}" Let $P$ be the largest $P\in\D$ as in (a). Then there exists $w\in\bold W(H)$ 
such that $P=S_w$.

\item "{(c)}" If $H$ is a $p$-subgroup of $\ca L$ then there exists $g\in\ca L$ such that 
$H^g\leq S$. 

\endroster 
\endproclaim 

\demo {Proof} For any $w=(h_1,\cdots,h_n)\in\bold W(H)$ let $w'$ be the word $(g_1,\cdots,g_n)$ defined by 
$g_i=h_1\cdots h_i$. As $H$ is finite we may choose $w$ so as to maximize the cardinality of the set 
$X=\{g_1,\cdots,g_n\}$. Supposing that $X\neq H$, let $g\in H-X$, and 
set $h=\Pi(w)\i g$. Then the set of entries of $(w\circ(h))'$ is $X\cup\{g\}$, contrary 
to the maximality of $X$. Thus $X=H$. 

We have $\bold W(\ca H)\sub\bold D$ as $H$ is a subgroup of $\ca L$, and thus $w\in\bold D$ via 
some $P\in\D$. Then $P^{g_i}=P^{h_1\cdots h_i}\leq S$ for all $i$, and so $P^h\leq S$ for 
all $h\in H$. Set $U=\<P^h\mid h\in H\>$ (the subgroup of $S$ generated by the union of all 
$P^h$ for $h\in H$). Then $U\in\D$ by (O2). In order to complete the proof of (a) and (b) it now suffices 
to show that $H\leq N_{\ca L}(U)$, and  for this it is enough to observe that, by 2.4(b), $(P^f)^g$ 
is defined and is equal to $P^{fg}$ for all $f,g\in H$. 

Next, by 2.10 there exists 
$V\in\D$ and $g\in\ca L$ such that $V=U^g$ and such that $N_S(V)\in Syl_p(N_{\ca L}(V))$. Let 
$c_g:N_{\ca L}(U)\to N_{\ca L}(V)$ be the isomorphism given by 2.3(b). Thus $H^g$ is a 
$p$-subgroup of $N_{\ca L}(V)$, so there exists $x\in N_{\ca L}(V)$ with $(H^g)^x\leq N_S(V)$. 
Since $(g,x)\in\bold D$ via $U$ we may apply 2.3(c), obtaining $(H^g)^x=H^{gx}$. Thus (c) holds 
with $gx$ in the role of $g$.  
\qed 
\enddemo

The theory being developed here purports to be ``$p$-local", so it is important to be able to analyze 
the structure of $N_{\ca L}(U)$ for $U$ an arbitrary subgroup of $S$ (not necessarily in $\D$). Since our 
policy in this Part I is to avoid bringing in the language of fusion 
systems, we shall obtain for now only the following very limited result concerning such normalizers.

\proclaim {Lemma 2.13} Let $(\ca L,\D,S)$ be a locality, let $R\leq S$ be a subgroup of $S$, and set 
$\G=\{N_P(R)\mid P\in\D\}$. Then:
\roster 

\item "{(a)}" $N_{\ca L}(R)$ is a partial subgroup of $\ca L$. 

\item "{(b)}" If $\G\sub\D$ then $(N_{\ca L}(R),\G)$ is an objective partial group. 

\item "{(c)}" If $\G\sub\D$ and $N_S(R)$ is a maximal $p$-subgroup of $N_{\ca L}(R)$ then 
$(N_{\ca L}(R),\G,S)$ is a locality. 

\endroster 
In particular, if $R\norm S$ then $(N_{\ca L}(R),\D,S)$ is a locality. 
\endproclaim 

\demo {Proof} Set $\ca L_R=N_{\ca L}(R)$. Then 2.3(c) shows that $\Pi$ maps $\bold W(\ca L_R)\cap\bold D$ 
into $\ca L_R$, while 1.6(c) shows that $\ca L_R$ is closed under the inversion in $\ca L$. Thus $\ca L_R$ 
is a partial subgroup of $\ca L$. 

Let $\bold D(\ca L_R)$ be the domain of the product in the partial group $\ca L_R$. Then 
$\bold D(\ca L_R)=\bold W(\ca L_R)\cap\bold D$, by definition. Assume now that $\G\sub\D$, and 
let $w\in\bold D(\ca L_R)$ via some $P\in\D$. Then $w\in\bold D(\ca L_R)$ via $N_P(R)$, and thus $\ca L_R$ 
satisfies the condition (O1) for objectivity in definition 2.1. The set $\G$ is evidently closed in the 
sense of condition (O2) in definition 2.1, so we have (b). Point (c) is then immediate from definition 2.8. 
\qed 
\enddemo

\proclaim {Lemma 2.14} Let $(\ca L,\D,S)$ be a locality, and set 
$$ 
O_p(\ca L)=\bigcap\{S_w\mid w\in\bold W(\ca L)\}. 
$$ 
Then $O_p(\ca L)$ is the unique largest subgroup of $S$ which is a partial normal subgroup of $\ca L$.  
\endproclaim 

\demo {Proof} Set $Y=O_p(\ca L)$, let $g\in\ca L$, and let $w\in\bold W(\ca L)$. Then $Y\leq S_g$ and 
$Y\leq S_{(g)\circ w}$, so $Y^g\leq S_w$. Thus $Y^g\leq Y$, and so $Y\norm\ca L$. 

Now let $X\leq S$ with $X\norm\ca L$. Then $X\norm S$. Let $g\in\ca L$, and set $P=S_g$. Then conjugation by 
$g$ is defined on $N_{\ca L}(P)$ by 2.3(b), so $N_X(P)^g$ is defined, and then $N_X(P)^g\leq X$ as 
$X\norm\ca L$. Thus $N_X(P)^g\leq S$, so $N_X(P)\leq S_g$. As $S_gX=PX$ is a subgroup of $S$, it follows 
that $X\leq P$, so $X^g=X$, and $X\leq O_p(\ca L)$. 
\qed 
\enddemo 

Given a locality $\ca L$, there can be more than one choice for $S$ and, after that, more than one choice 
for $\D$, such that $\bold D=\bold D_\D$.

\definition {Definition 2.15} Let $\ca L$ be a partial group. An automorphism $\a$ of $\ca L$ is {\it inner}  
if there exists $g\in\ca L$ such that $\a$ is given by conjugation by $g$. That is: 
\roster 

\item "{(1)}" $x\in\ca L$ $\implies$ $x^g$ is defined and is equal to $x\a$, and  

\item "{(2)}" $(x_1,\cdots,x_n)\in\bold D$ $\implies$ $(x_1^g,\cdots,x_n^g)\in\bold D$ 
and $x_1^g\cdots x_n^g=(x_1\cdots x_n)^g$. 

\endroster 
Write $Inn(\ca L)$ for the set of all inner automorphisms of $\ca L$. 
\enddefinition

\definition {Definition 2.16} Let $\ca L$ be a locality and let $S$ be a $p$-subgroup of $\ca L$. 
Then $S$ is a {\it Sylow $p$-subgroup} of $\ca L$ if there exists a set $\D$ of subgroups of $S$ 
such that $(\ca L,\D,S)$ satisfies the conditions (L1) and (L2) of definition 2.8. Write 
$Syl_p(\ca L)$ for the set of Sylow $p$-subgroups of $\ca L$. 
\enddefinition 

It will be convenient to introduce some notation regarding ``conjugation" by words $w\in\bold(\ca L)$, 
in the case of a locality $(\ca L,\D,S)$. Thus, let $w=(g_1,\cdots,g_n)\in\bold W(\ca L)$, 
and let $X$ be a subgroup of $S_w$. For each $x=x_0\in X$ and each index $i$ from $1$ to $n$ there is 
then an element $x_i\in S$ defined by $x_i=(x_{i-1})^{g_i}$, and we may write $x^w$ for $x_n$. 
The mapping $x\maps x^w$ may be written as
$$ 
c_w: X\to X^w. 
$$

\proclaim {Proposition 2.17} Let $(\ca L,\D,S)$ be a locality, and let $K$ be 
the set of all $g\in\ca L$ such that conjugation by $g$ is an inner automorphism of $\ca L$. 
\roster 

\item "{(a)}" $K$ is a subgroup of $\ca L$.

\item "{(b)}" For each $w=(x_1,\cdots,x_n)\in\bold D$ and each $(g_0,\cdots,g_n)\in\bold W(K)$, 
we have 
$$ 
(g_0,x_1,g_1\i,g_1,x_2,\cdots,g_{n-1}\i,g_{n-1},x_n,g_n)\in\bold D. \tag*
$$

\item "{(c)}" Set $Q=O_p(K)$. Then $Q\in\D$ and $K=N_{\ca L}(Q)$. Moreover, for each $w\in\bold D$ 
there exists $P\leq Q\cap S_w$ such that $P\in\D$ and $P^w\leq Q$. 

\item "{(d)}" $Syl_p(\ca L)=Syl_p(K)$. 

\item "{(e)}" $Inn(\ca L)\norm Aut(\ca L)$, and every automorphism of $\ca L$ can be factored as an 
inner automorphism followed by an automorphism which leaves $S$ invariant. 

\endroster 
\endproclaim 

\demo {Proof} Write $\ca S$ for $Syl_p(\ca L)$. For $S'\in\ca S$ note that $N_{\ca L}(S')$ is a 
subgroup of $\ca L$. By 2.9 the condition (*) holds with $N_{\ca L}(S')$ in the role of $K$. 
Let $\ca U$ be the union (taken over all $S'\in\ca S$) of the groups $N_{\ca L}(S')$, and let $H$ be the 
partial subgroup of $\ca L$ generated by $\ca U$. The above observation concerning 2.9 together with a 
straightforward argument by induction on word-length then yields $\bold W(\ca U)\sub\bold D$, 
and thus $H$ is a subgroup of $\ca L$. Moreover, $\bold W(\ca U)$ has the following property:  
For any $(g_1,\cdots,g_n)\in\bold D(\ca L)$, and any $(u_0,\cdots,u_n)$ with each $u_i\in\bold W(\ca U)$,  
the word  
$$ 
u_0\circ g_1\circ u_1\i\circ u_1\circ g_2\circ\cdots\circ u_{n-1}\i\circ u_{n-1}\circ g_n\circ u_n 
$$ 
is in $\bold D$, by induction on the sum of the lengths of the words $u_i$. The condition 
(*) then holds with $H$ in the role of $K$, by $\bold D$-associativity. In particular, conjugation 
by $g\in H$ is an automorphism of $\ca L$, and thus $H$ is a subset of $K$. 

Now let $g\in K$. Then $(\ca L,\D^g,S^g)$ is a locality, so $S^g\in\ca S$. Then $S^g\in Syl_p(H)$, and 
there exists $h\in H$ with $S^g=S^h$. The product $gh\i$ is defined, and 
then $gh\i\in N_{\ca L}(S)$. Thus $gh\i\in H$, so $g\in H$, and we conclude that $H=K$. This completes 
the proof of (a), (b), and (d). Point (e) is immediate from (d), so it remains only to prove (c).   

By 2.11(a), $K$ normalizes a member of $\D$. Then, since $S\in Syl_p(K)$, we obtain $Q:=O_p(K)\in\D$. 
By 2.11(b) there exists a word $u\in\bold W(K)$ such that 
$Q=S_u$. Then also $Q=S_{u\i}$. Let $w\in\bold D$, and let $f\in N_{\ca L}(Q)$. 
The word $v:=u\circ w\circ u$ is in $\bold D$ by (*), and  
clearly $S_v\leq Q$ and $S_{v\i}\leq Q$. Set $P=(S_v)^{\Pi(u)}$. Then $P\leq S_w\cap Q$ and 
$P^{\Pi(w)}\leq Q$. As $S_v\in\D$ we have $P\in\D$. Write $w=(x_1,\cdots,x_n)$. Then 
$$ 
(f\i,x_1,f,f\i,\cdots,f,f\i,x_n,f)\in\bold D\ \ \text{via $P$}, 
$$ 
and this shows that conjugation by $f$ is an inner automorphism of $\ca L$. Thus $K=N_{\ca L}(Q)$, 
and the proof is complete.  
\qed 
\enddemo

\proclaim {Lemma 2.18} Let $\ca L$ be a locality and let $S$ be a Sylow $p$-subgroup of $\ca L$. 
Then there is a unique smallest set $\D=\D_0$ and a unique largest set $\D=\D_1$ of subgroups of $S$ 
such that the conditions (L1) and (L2) of 2.8 are satisfied by $(\ca L,\D,S)$.  
\endproclaim 

\demo {Proof} Set Take $\D_0$ to be the overgroup-closure in $S$ of the set of all $S_w$ for 
$w\in\bold D$. Take $\D_1$ to be the union of all the the sets $\G$ of subgroups of $S$ which fulfill 
(L1) and (L2). 
\qed 
\enddemo

\vskip .2in 
\noindent 
{\bf Section 3: Partial normal subgroups} 
\vskip .1in 

Throughout this section we fix a locality $(\ca L,\D,S)$ and a partial normal subgroup $\ca N\norm\ca L$. 
Recall that this means that $\ca N\leq\ca L$ is a partial subgroup of $\ca L$ and that 
$\Pi(g\i,x,g)\in\ca N$ for all $x\in\ca N$ and all $g\in\ca L$ for which $(g\i,x,g)\in\bold D$.  
Set $T=S\cap\ca N$.  

\proclaim {Lemma 3.1} 
\roster

\item "{(a)}" $(T\cap S_g)^g\leq T$ for all $g\in\ca L$. 

\item "{(b)}" Let $x\in\ca N$ and let $P$ be a subgroup of $S_x$. Then $PT=P^x T$.

\item "{(c)}" $T$ is maximal in the poset of $p$-subgroups of $\ca N$. 

\endroster
\endproclaim

\demo {Proof} (a) Let $g\in\ca L$ and let $t\in S_g\cap T$. Then $t^g\in S$, and $t^g\in\ca N$ as 
$\ca N\norm\ca L$. Thus $t^g\in T$. 
\vskip .1in 
\noindent 
(b) Let $a\in P$. Then $(P^x)^a\leq S$ and $P^a=P$. Setting $w=(a\i,x\i,a,x)$ we 
then have $w\in\bold D$ via $P^{xa}$. Now $\Pi(w)=a\i a^x\in S$, while also 
$\Pi(w)=(x\i)^ax\in\ca N$, and so $\Pi(w)\in T$. Then $a^x\in aT$, and we have thus 
shown that $P^x\leq PT$. Then $P^x T\leq PT$. The equality $P^xT=PT$ can then be 
deduced from (a) (which implies that $|P^x\cap T|=|P\cap T|$), or from symmetry with 
$x\i$ and $P^x$ in place of $x$ and $P$. 
\vskip .1in 
\noindent 
(c) Let $R$ be a $p$-subgroup of $\ca N$ containing $T$. By 2.9(b) there exists 
$g\in\ca L$ with $R^g\leq S$, and then $R^g\leq S\cap\ca N=T$. As $T^g=T$ and conjugation by $g$ is 
injective, we conclude that $R=T$. 
\qed
\enddemo

\proclaim {Lemma 3.2} Let $x,y\in\ca N$ and let $f\in N_{\ca L}(T)$. 
\roster 

\item "{(a)}" If $(x,f)\in\bold D$ then $(f,f\i,x,f)\in\bold D$, $xf=fx^f$, 
and $S_{(x,f)}=S_{(f,x^f)}=S_x\cap S_f$. 

\item "{(b)}" If $(f,y)\in\bold D$ then $(f,y,f\i,f)\in\bold D$, 
$fy=y^{f\i}f$, and $S_{(f,y)}=S_{(y^{f\i},f)}=S_{y^{f\i}}\cap S_f$. 

\endroster 
\endproclaim 

\demo {Proof} Set $Q=S_{(x,f)}$ and $P=S_x\cap S_f$. As $x\in\ca N$ we have $Q^xT=QT$ by 3.1(b).  
Then since $f\in N_{\ca L}(T)$ we obtain $QT\leq S_f$. Thus $Q\leq P$. But also $P^xT=PT$, so $P=Q$. 
Now $(f,f\i,x,f)\in\bold D$ via $Q$, and $\Pi(f,f\i,x,f)=xf=fx^f$. As  
$$ 
Q=S_f\cap S_x = S_{(f,f\i,x,f)}\leq S_{(f,x^f)}, 
$$ 
we obtain $S_{(x,f)}\leq S_{(f,x^f)}$. Thus, in order to complete the proof of (a), it remains to show 
that $S_{(f,x^f)}\leq Q$. 

On the other hand, in addressing (b), set $R=S_{(f,y)}$. Then $R^{fy}T=R^fT\leq S_{f\i}$, so 
$(f,y,f\i,f)\in\bold D$ via $R$, and hence $R=S_{(f,y,f\i,f)}$. Take $y=x^f$. Thus 
$$ 
R=S_{(f,x^f)}=S_{(f,x^f,f\i,f)}\leq S_{(x,f)}=Q,  
$$ 
completing the proof of (a). Moreover, we have shown that $fy=y^{f\i}f$. The remainder of (b) now 
follows as an application of (a) to $(y^{f\i},f)$. The remainder of (b) now 
follows as an application of (a) to $(y^{f\i},f)$. 
\qed 
\enddemo

\proclaim {Lemma 3.3} Let $w\in\bold W(N_{\ca L}(T))$, set $g=\Pi(w)$, 
and let $x,y\in\ca N$. 
\roster 

\item "{(a)}" Suppose that $(x)\circ w\in\bold D$ and set $P=S_{(x)\circ w}$. 
Then $u:=w\i\circ(x)\circ w\in\bold D$, and $S_u=P^g$. 

\item "{(b)}" Suppose that $w\circ(y)\in\bold D$ and set $Q=S_{w\circ(y)}$. Then 
$v:=w\circ(y)\circ w\i\in\bold D$, and $S_v=Q$. 

\endroster 
\endproclaim 

\demo {Proof} We prove only (a), leaving it to the reader to supply a similar argument for (b). 
As $P^xT=PT$ by 3.1(b), and since both $P^x$ and $T$ are contained in $S_w$, we have $P\leq S_w$. 
Then $P^g\leq S$ by 2.3(c), and $P^g\leq S_u$. In particular $u\in\bold D$ via $P^g$. As  
$S_u\leq S_{w\i}$, and $(S_u)^{g\i}\leq P$, we obtain $S_u\leq P^g$. Thus $S_u=P^g$.  
\qed 
\enddemo

\proclaim {Lemma 3.4} Let $w=(f_1,g_1,\cdots,f_n,g_n)\in\bold D$, with 
$f_i\in N_{\ca L}(T)$ and with $g_i\in\ca N$ for all $i$. Set $u=(f_1,\cdots,f_n)$. Then there exists 
$g\in\ca N$ such that 
$$ 
S_w\leq S_{u\circ(g)}\ \ \text{and}\ \ \Pi(w)=\Pi(u\circ(g)). \tag*
$$  
Similarly, there exists $g'\in\bold{\ca N}$ such that 
$$ 
S_w\leq S_{(g')\circ u}\ \ \text{and}\ \ \Pi(w)=\Pi((g')\circ u). \tag**
$$ 
\endproclaim 

\demo {Proof} The case $n=1$ is given by 3.2. For the general case, write 
$$ 
w=(f_1,g_1)\circ w_1,  
$$  
and set $u_1=(f_2,\cdots,f_n)$. Induction on $n$ implies that there exists $x_1\in\ca N$ such that 
$S_{w_1}\leq S_{u_1\circ(x_1)}$ and $\Pi(w_1)=\Pi(u_1\circ(x_1))$. Set 
$$ 
w'=(f_1,g_1)\circ u_1\circ(x_1). 
$$ 
Thus $S_w\leq S_{w'}$, and we have $\Pi(w)=\Pi(w')$ by $\bold D$-associativity. Apply 3.3(a) to the word 
$(g_1)\circ u_1$ to obtain (*). A similar argument in which one begins by writing 
$w=w_n\circ(f_n,g_n)$ yields (**).  
\qed 
\enddemo 

The following result will be of fundamental importance in Part III. 

\proclaim {Lemma 3.5} Suppose that 
$$ 
C_{N_{\ca L}(P)}(O_p(N_{\ca L}(P))\leq O_p(N_{\ca L}(P)\quad\forall\  P\in\D. \tag* 
$$ 
Then $N_{\ca N}(T)\leq N_{\ca L}(C_S(T)T)\leq N_{\ca L}(C_S(T))$. 
\endproclaim 

\demo {Proof} Let $g\in N_{\ca N}(T)$ and set $P=S_g$. Then $P^g=P$ by 3.1(b). Set $M=N_{\ca L}(P)$, 
$K=\ca N\cap M$, and $D=N_{C_S(T)}(P)$. Here $M$ is a subgroup of $\ca L$ as $P\in\D$, and then 
$K$ is a normal subgroup of $M$ by 1.7(d). We have $T\norm M$ by 3.1(a), so $D\leq C_M(T)$, and then 
$[K,D]\leq C_K(T)$. Notice that $T\in Syl_p(K)$ by 3.1(c), so that $K/T$ is a $p'$-group. Then 
$C_K(T)=Z(T)\times Y$ where $Y$ is a $p'$-group; and thus $Y$ is a normal $p'$-subgroup of $M$. Then 
$[O_p(M),Y]=1$, so $Y\leq O_p(M)$ by hypothesis, and thus $Y=1$. Thus $[K,D]\leq T$, so $[g,D]\leq T$, and 
$D\leq P$. As $C_S(T)P$ is a $p$-group we thereby obtain $C_S(T)\leq P$. This shows that  
$C_S(T)$ is $g$-invariant for all $g\in\ca N_{\ca N}(T)$, and yields the lemma.   
\qed 
\enddemo

\definition {Definition 3.6} Let $\ca L\circ\D$ be the set of all pairs 
$(f,P)\in\ca L\times\D$ such that $P\leq S_f$. Define a relation $\uparrow$  
on $\ca L\circ\D$ by $(f,P)\uparrow(g,Q)$ if there exist elements 
$x\in N_{\ca N}(P,Q)$ and $y\in N_{\ca N}(P^f,Q^g)$ such that $xg=fy$. 
\enddefinition 

This relation may be indicated by means of a commutative square: 
$$ 
\CD 
Q@>g>>Q^g \\ 
@AxAA  @AAyA  \\ 
P@>f>>P^f 
\endCD \tag* 
$$
of conjugation maps, labeled by the conjugating elements, and in which the 
horizontal arrows are isomorphisms and the vertical arrows are injective 
homomorphisms. The relation $(f,P)\up(g,Q)$ may also be expressed by: 
$$
\text{$w:=(x,g,y\i,f\i)\in\bold D$ via $P$, and $\Pi(w)=\1$.} 
$$ 
\vskip .1in
It is easy to see that $\up$ is reflexive and transitive. We say that $(f,P)$  
is {\it maximal} in $\ca L\circ\D$ if $(f,P)\up(g,Q)$ implies that $|P|=|Q|$.  
As $S$ is finite there exist maximal elements in $\ca L\circ\D$. Since 
$(f,P)\up(f,S_f)$ for $(f,P)\in\ca L\circ\D$, we have $P=S_f$ for every 
maximal $(f,P)$. For this reason, we will say that $f$ is {\it $\up$-maximal}
in $\ca L$ (with respect to $\ca N$) if $(f,S_f)$ is maximal in 
$\ca L\circ\D$.

\proclaim {Lemma 3.7} Let $f\in\ca L$. 
\roster 

\item "{(a)}" If $f\in N_{\ca L}(S)$ then $f$ is $\up$-maximal. 

\item "{(b)}" If $f$ is $\up$-maximal then so is $f\i$. 

\item "{(c)}" If $f$ is $\up$-maximal and $(f,S_f)\up(g,Q)$, then $g$ is $\up$-maximal and 
$Q=S_g$. 

\endroster  
\endproclaim 

\demo {Proof} Point (a) is immediate from definition 3.6. Now suppose that $f$ is $\up$-maximal, 
and let $g\in\ca L$ with $(f\i,S_{f\i})\up(g\i,S_{g\i})$. Since $S_{f\i}=(S_f)^f$ and 
$S_{g\i}=(S_g)^g$, one obtains a diagram 
$$ 
\CD 
(S_g)^g@>g\i>>S_g \\ 
@AxAA      @AAyA \\ 
(S_f)^f@>f\i>>S_f 
\endCD 
$$ 
as in definition 3.6, from which it is easy to read off the relation 
$(f,S_f)\up(g,S_g)$. Then $|S_f|=|S_g|$ as $f$ is $\up$-maximal, and then 
also $|S_{f\i}|=|S_{g\i}|$. Thus $f\i$ is $\up$-maximal, and (b) holds. Point (c) is 
immediate from the transitivity of $\up$. 
\qed 
\enddemo 

\proclaim {Lemma 3.8} Let $(g,Q),(h,R)\in\ca L\circ\D$ with$(g,Q)\up(h,R)$, and suppose that 
$T\leq R$. Then there exists a unique $y\in\ca N$ with $g=yh$. Moreover: 
\roster 

\item "{(a)}" $Q^y\leq R$, and $Q\leq S_{(y,h)}$. 

\item "{(b)}" If $N_T(Q^g)\in Syl_p(N_{\ca N}(Q^g))$, then $N_T(Q^y)\in Syl_p(N_{\ca N}(Q^y))$. 

\endroster
\endproclaim

\demo {Proof} By definition of the relation $\up$, there exist elements $u\in N_{\ca N}(Q,R)$ and 
$v\in N_{\ca N}(Q^g,R^h)$ such that $(u,h,v\i,g\i)\in\bold D$ via $Q$, and such that $\Pi(w)=\1$. 
$$ 
\CD 
R@>h>>R^h \\ 
@AuAA  @AAvA \\  
Q@>>g>Q^g 
\endCD 
$$
In particular, $uh=gv$. Since $T\leq R$, points (a) and (b) of 3.1 yield 
$$
\text{$T=T^h$, $Q^uT=QT\leq R$, and $Q^gT=Q^{gv}T\leq R^h$}. 
$$
Then
$$
w:=(u,h,v\i,h\i)\in\bold D\quad{\text via}\quad (Q,Q^u,Q^{uh},Q^{uhv\i}=Q^g,Q^{gh\i}). 
$$
Set $y=\Pi(w)$. Then $y=u(v\i)^{h\i}\in N_{\ca N}(Q,R)$. Since 
$(u,h,v\i,h\i,h)$ and $(g,v,v\i)$ are in $\bold D$ (as $\ca L$ is a partial 
group), we get $yh=uhv\i=g$. This yields (a). The uniqueness of $y$ is 
given by right cancellation. 

Suppose now that $N_T(Q^g)\in Syl_p(N_{\ca N}(Q^g))$. As 
$N_T(Q^y)^h=N_T(Q^g)$, it follows from 2.3(b) that $N_T(Q^y)\in Syl_p(N_{\ca N}(Q^y))$. 
\qed 
\enddemo

\proclaim {Proposition 3.9} Let $g\in\ca L$ and suppose that $g$ is $\up$-maximal with respect to 
$\ca N$. Then $T\leq S_g$. 
\endproclaim 

\demo {Proof} Set $P=S_g$ and $Q=P^g$. We first show: 
\roster 

\item "{(1)}" Let $y\in N_{\ca N}(P,S)$. Then $|T\cap P|=|T\cap P^y|$, and $(g,P)\up(y\i g,P^y)$. 
In particular, $y\i g$ is $\up$-maximal.

\endroster 
As $P^y T=PT$ we obtain 
$$ 
|P^y:P^y\cap T|=|P^yT:T|=|PT:T|=|P:P\cap T|, 
$$ 
and so $|T\cap P|=|T\cap P^y|$. The following diagram 
$$ 
\CD 
P^y @>y\i g>> P^g \\ 
@AyAA     @AA\1A  \\ 
P  @>>g> P^g 
\endCD \tag*
$$ 
shows that $(g,P)\up(y\i g,P^y)$ and completes the proof of (1). 

Suppose next that $N_T(P)\in Syl_p(N_{\ca N}(P))$. Then $N_T(P)^g\in Syl_p(N_{\ca N}(Q))$, by 2.3(b), and 
there exists $x\in N_{\ca N}(Q)$ such that $N_T(Q)^x\leq N_T(P)^g$. Here $(x,g\i)\in\bold D$ via $Q$, and  
we get $N_T(Q)^{xg\i}\leq N_T(P)\leq S$. Thus $(g\i,Q)\up(xg\i,N_T(Q)Q)$. As $g\i$ is $\up$-maximal by 3.7, 
it follows that $N_T(Q)\leq Q$. Hence $T\leq Q$ and then $T\leq P$. We have thus shown:  
\roster 

\item "{(2)}" If $T\nleq P$ then $N_T(P)\notin Syl_p(N_{\ca N}(P))$. 

\endroster 

We next show: 
\roster 

\item "{(3)}" Suppose that there exists $y\in\ca N$ such that $P\leq S_y$ and such that 
$N_T(P^y)\in Syl_p(N_{\ca N}(P^y))$. Then $T\leq P$. 

\endroster 
Indeed, under the hypothesis of (3) we have $(g,P)\up(y\i g,P^y)$ by (1). Then 
$(y\i g,P^y)$ is $\up$-maximal and $P^y=S_{y\i g}$. If $T\nleq P$ then $T\nleq P^y$, and then (2) 
applies to $(y\i g,P^y)$ in the role of $(g,P)$ and yields a contradiction. So, (3) holds. 

Among all counter-examples, let $g$ be chosen so that  $|P|$ is as large as possible. By 2.10 there 
exists $f\in\ca L$ so that $Q^f\leq S$ and so that $N_S(Q^f)\in Syl_p(N_{\ca L}(Q^f))$. Set $h=gf$ 
(where the product is defined via $P$) and set $R=P^h$. Let $(h',P')$ be maximal in $\ca L\circ\D$ 
with $(h,P)\up(h',P')$, and set $R'=R^{h'}$. If $T\leq P'$ then 3.8 yields $h=yh'$ for some $y\in\ca N$ 
such that $P\leq S_y$ and such that $N_T(P^y)\in Syl_p(N_{\ca N}(P^y))$. The existence of 
such an element $y$ contradicts (3), so we conclude that $T\nleq P'$. Then $(h',P')$ is a counter-example 
to the proposition, and the maximality of $|P|$ yields $|P|=|P'|$. Then $(h,P)$ is maximal in 
$\ca L\circ\D$, as is $(h\i,R)$ by 3.7(b), and so $T\nleq R$. But $N_T(R)\in Syl_p(N_{\ca N}(R))$ since 
$R=Q^f$, so (2) applies with $(h\i,R)$ in the role of $(g,P)$, and yields $T\leq R$. Then $T\leq P$, 
and the proof is complete. 
\qed
\enddemo

\proclaim {Lemma 3.10} Suppose that $N_{\ca N}(T)\leq N_{\ca L}(S)$. 
Then every element of $N_{\ca L}(T)$ is $\up$-maximal with respect to $\ca N$. 
\endproclaim 

\demo {Proof} Let $f\in N_{\ca L}(T)$ and set $P=S_f$. Then $T\leq P$ and $T\leq P^f$. Let 
$(g,Q)\in\ca L\circ\D$ with $(f,P)\up(g,Q)$, and let $x,y\in\ca N$ be chosen as in definition 3.6. Then 
$P\leq Q$ by 3.9, and we have $P=P^x$ and $(P^f)^y=P^f\leq Q^g$ by 3.1(b). In order to show that $f$ is 
$\up$-maximal it suffices now to show that $P=Q$, and hence it suffices to show that $N_Q(P)\leq P$. 

Set $D=N_Q(P)$. As $x\in N_{\ca N}(T)\leq N_{\ca L}(S)$ we have $[D,x]\leq S\cap\ca N=T$, and similarly 
$[D^g,y]\leq T$. This shows that $(x,g,y\i)\in\bold D$ via $D$. As $xg=fy$ 
by the setup of definition 3.6, we have $xgy\i=f$ by 1.3(g), and thus $D\leq S_f$. That is, we have 
$N_Q(P)\leq P$, as required. 
\qed 
\enddemo

If $X$ and $Y$ are subsets of $\ca L$ then one has the notion of the product $XY$, introduced in section 1, 
as the set of all $\Pi(x,y)$ with $(x,y)\in\bold D\cap(X\times Y)$.

\proclaim {Corollary 3.11 (Frattini Lemma)} Let $(\ca L,\D,S)$ be a locality, let $\ca N\norm\ca L$ be a 
partial normal subgroup, and let $\L$ be the set of $\up$-maximal elements of $\ca L$ with 
respect to $\ca N$. Then $\ca L=\ca N\L=\L\ca N$. In particular, we have 
$\ca L=N_{\ca L}(T)\ca N=\ca N N_{\ca L}(T)$. 
\endproclaim 

\demo {Proof} Let $f\in\ca L$, set $P=S_f$, and choose $(g,Q)\in\ca L\circ\D$ so that 
$(f,P)\up(g,Q)$ and so that $g$ is $\up$-maximal. Then 3.8 yields $f=xg$ for some $x\in\ca N$, 
and then 3.2 shows that $f=gy$ where $y=x^g$. 
\qed 
\enddemo

The following result will be seen to play a crucial role in the theory being developed here. It was 
discovered and proved by Bernd Stellmacher, in his reading of an early draft of [Ch1]. The proof given 
here is his.

\proclaim {Lemma 3.12 (Splitting Lemma)} Let $(x,f)\in\bold D$ with $x\in\ca N$ and with 
$f$ $\up$-maximal with respect to $\ca N$. Then $S_{(x,f)}=S_{xf}=S_{(f,x^f)}$. 
\endproclaim 

\demo {Proof} Appealing to 3.2: Set $y=x^f$ and $g=xf$ (so that also $g=fy$), and set $Q=S_{(x,f))}$ (so 
that also $Q=S_{(f,y)}$). Thus $Q\leq S_f\cap S_g$. Also, 3.2(a) yields $Q=S_f\cap S_x$.  Set 
$$ 
P_0=N_{S_f}(Q),\ \ P_1=N_{S_g}(Q),\ \ P=\<P_0,P_1\>, 
$$ 
and set $R=P_0\cap P_1$. Then $Q\leq R$. In fact, 2.3(b) shows that $y=f\i g$ and that 
$(R^f)^y=R^g$, so $R\leq Q$, and thus $P_0\cap P_1=Q$. Assume now that $(x,f)$ is a 
counter-example to the lemma. That is, assume $Q<S_g$ (proper inclusion). Then 
$Q<P_1$ and so $P_1\nleq P_0$. Thus: 
\roster 

\item "{(1)}" $P_1\nleq S_f$. 

\endroster 
Among all counter-examples, take $(x,f)$ so that $|Q|$ is as large as possible. We consider two 
cases, as follows. 
\vskip .1in 
\noindent 
CASE 1: $x\in N_{\ca L}(T)$. 
\vskip .1in 
As $f\in N_{\ca L}(T)$ (3.5) we have $T\leq Q$, and then $x\in N_{\ca L}(Q)$ by 3.1(b). 
Thus $Q^g=Q^{xf}=Q^f$. Set $Q'=Q^g$. Then 2.2(b) yields an isomorphism 
$c_f:N_{\ca L}(Q)\to N_{\ca L}(Q')$. Here $f=x\i g$ so $c_f=c_{x\i}\circ c_g$ by 2.2(c). 
As $x\in N_{\ca N}(Q)\norm N_{\ca L}(Q)$, we obtain $(P_1)^{x\i}\leq N_{\ca N}(Q)P_1$, and then 
$$ 
(P_1)^f=((P_1)^{x\i})^g\leq (N_{\ca N}(Q)P_1)^g\leq N_{\ca N}(Q')N_S(Q'). 
$$ 
Also $(P_0)^f\leq N_S(Q')$, so 
\roster 

\item "{(2)}" $P^f\leq N_{\ca N}(Q')N_S(Q')$. 

\endroster  
Since $T\leq Q'$, $T$ is a Sylow $p$-subgroup of $N_{\ca N}(Q')$ by 3.1(c), and thus 
$N_S(Q')$ is a Sylow $p$-subgroup of $N_{\ca N}(Q')N_S(Q')$. By (2) and Sylow's 
Theorem there is then an element $v\in N_{\ca N}(Q')$ such that $P^{fv}\leq N_S(Q')$. In 
particular, we have: 
\roster 

\item "{(3)}" $P_0\leq S_{(f,v)}$ and $P\leq S_{fv}$. 

\endroster 
Set $u=v^{f\i}$. Then $(u,f)\in\bold D$ and we have $uf=fv$ and $S_{(u,f)}=S_{(f,v)}$ by 3.2. If  
$S_{(f,v)}=S_{fv}$ then $(f,v)\in\bold D$ via $P$, so that $P\leq S_f$, contrary to (1). Thus 
$S_{(f,v)}\neq S_{fv}$, and so $(u,f)$ is a counter-example to the lemma. Then (3) and the maximality of 
$|Q|$ in the choice of $(x,f)$ yields $Q=P_0=N_{S_f}(Q)$, and so $Q=S_f$. As 
$(f,Q)\up(g,P)$ via $(x\i,\1)$, we have contradicted the $\up$-maximality of $f$. 
\vskip .1in 
\noindent 
CASE 2: The case $x\notin N_{\ca L}(T)$. 
\vskip .1in 
Let $h$ be $\up$-maximal, with $(g,S_g)\up(h,S_h)$. Then $T\leq S_h$ by 3.9, and there exists  
$r\in\ca N$ with $g=rh$ by 3.8. Note that 3.8 yields also $S_g=S_{(r,h)}\geq Q=S_{(x,g)}$.  

Set $w=(f\i,x\i,r,h)$, observe that $w\in\bold D$ via $Q^g$, and find 
$$ 
\Pi(w)=(f\i x\i)(rh)=g\i g=\1. 
$$ 
Then 2.3 yields $h=r\i xf$. Since both $f$ and $h$ are in 
$N_{\ca L}(T)$, 2.2(c) yields $r\i x\in N_{\ca L}(T)$, and so $r\i x\in N_{\ca N}(T)$. Then  
Case 1 applies to $(r\i x,f)$, and thus $S_h=S_{(r\i x,f)}\leq S_f$ (using 3.2). By definition of 
$\up$ there exist $a,b\in\ca N$ such that one has the usual  ``commutative diagram": 
$$ 
\CD 
S_h@>h>>S_{h\i} \\ 
@AaAA   @AAbA  \\ 
S_g@>>g>S_{g\i} 
\endCD. 
$$ 
As $T\leq S_h$, 3.1(b) yields 
$$ 
S_g\leq S_gT=(S_g)^aT\leq S_h, 
$$ 
and so $S_g\leq S_f$. This again contradicts (1), and completes the proof. 
\qed
\enddemo

The splitting lemma yields a useful criterion for partial normality, as follows.

\proclaim {Corollary 3.13} Let $\ca L$ be a locality, let 
$\ca N\norm\ca L$, and let $\ca K\norm\ca N$ be a partial normal subgroup of 
$\ca N$. Suppose that $\ca K$ is $N_{\ca L}(T)$-invariant. I.e. suppose that 
$x^h\in\ca K$ for all $(h\i,x,h)\in\bold D$ such that $x\in\ca K$ and 
$h\in N_{\ca L}(T).)$ Then $\ca K\norm\ca L$. 
\endproclaim 

\demo {Proof} Let $x\in\ca K$ and let $f\in\ca L$ such that $x^f$ is defined. 
By the Frattini Lemma we may write $f=yg$ with $y\in\ca N$ and with $g$ $\up$-maximal, 
and then the splitting lemma yields $S_f=S_{(y,g)}$. Set $u=(f\i,x,f)$ and 
$v=(g\i,y\i,x,y,g)$. Then $S_u=S_v\in\D$, and $x^f=\Pi(u)=\Pi(v)=(x^y)^g$. Thus 
$x^f\in\ca K$, and $\ca K\norm\ca L$. 
\qed 
\enddemo

A subset $X$ of $\ca L$ of the form $\ca N f$, $f\in\ca L$, will be called a {\it coset} of $\ca N$. 
A coset $\ca Nf$ is {\it maximal} if it not a proper subset of any coset of $\ca N$.

\proclaim {Proposition 3.14} The following hold. 
\roster 

\item "{(a)}" $\ca Nf=f\ca N$ for all $f\in N_{\ca L}(T)$, and if $f$ is $\up$-maximal with respect to 
$\ca N$ then $\ca Nf=\ca N f\ca N=f\ca N$. 

\item "{(b)}" Let $f,g\in\ca L$ and assume that $f$ is $\up$-maximal. Then  
$$ 
(g,S_g)\up(f,S_f)\ \iff\ \ca Ng\sub\ca Nf\ \iff\ g\in\ca Ng.   
$$ 

\item "{(c)}" If $g\in\ca L$ is $\up$-maximal relative to $\ca N$ then $\ca Ng$ is a maximal coset 
of $\ca N$. Moreover, every maximal coset of $\ca N$ is of the form $\ca Nf$ for some $\up$-maximal 
element $f$.  

\item "{(d)}" $\ca L$ is partitioned by the set $\ca L/\ca N$ of maximal cosets of $\ca N$. 

\item "{(e)}" Let $u:=(g_1,\cdots,g_n)\in\bold D$ and let $v:=(f_1,\cdots,f_n)$ be a sequence of 
$\up$-maximal elements of $\ca L$ such that $g_i\in\ca Nf_i$ for all $i$. Then 
$TS_u\leq S_v$, and $\Pi(u)\in\ca N\Pi(v)$.  

\endroster 
\endproclaim 

\demo {Proof} (a): That $\ca Nf=f\ca N$ for $f\in N_{\ca L}(T)$ is given by 3.2. Now let $f$ be 
$\up$-maximal relative to $\ca N$. Then $f\in N_{\ca L}(T)$ by 3.9. Let $x,y\in\ca N$ such that 
$(x,f,y)\in\bold D$. Then $(x,fy)\in\bold D$, and $(x,fy)=(x,y'f)$ where $y'=fyf\i$. The splitting 
lemma (3.12) then yields $(x,y',f)\in\bold D$, and thus $\ca Nf\ca N\sub\ca Nf$. The reverse inclusion 
is obvious, and yields (a). 

\vskip .1in 
\noindent 
(b): We have $g=xf$ for some $x\in\ca N$ by 3.11, 
and then 3.12 yields $S_g=S_{(x,f)}\leq S_f$. Let $y\in\ca N$ such that $(y,g)\in\bold D$. 
Since $(S_{(y,g)})^y\leq S_g=S_{(x,f)}$ we get $(y,x,f)\in\bold D$, and $yg=(yx)f\in\ca Nf$. Thus 
$(g,S_g)\up(f,S_f)\implies \ca Ng\sub\ca Nf$. Clearly $\ca Ng\sub\ca Nf\implies g\in\ca Nf$. The 
required circle of implications is then completed by 3.12. 
\vskip .1in 
\noindent 
(c): Immediate from (b).  
\vskip .1in 
\noindent 
(d): Let $f$ and $g$ be $\up$-maximal, and let $h\in\ca Nf\cap\ca Ng$. Thus there exist $x,y\in\ca N$ 
with $(x,f)\in\bold D$, $(y,g)\in\bold D$, and with $h=xf=yg$. Then also $(x\i,x,f)\in\bold D$, so 
$(x\i,yg)=(x\i,xf)\in\bold D$. The splitting lemma (3.12) then yields $(x\i,y,g)\in\bold D$, and 
we thereby obtain $f=x\i yg\in\ca Ng$. Now (b) implies that $\ca Nf\sub\ca Ng$, and symmetry gives 
the reverse inclusion. Thus $\ca Nf=\ca Ng$ if $\ca Nf\cap\ca Ng\neq\nset$. 
\vskip .1in
\noindent
(e): Let $x_i\in\ca N$ with $g_i=x_if_i$. Set $Q_0=S_u$, and set 
$Q_i=(Q_{i-1})^{g_i}$ for $1\leq i\leq n$. Then 
$Q_{i-1}\leq S_{g_i}=S_{(x_i,f_i)}\leq S_{f_i}$ by 3.12, and then also 
$Q_{i-1}T=(Q_{i-1})^{x_i}T\leq S_{f_i}$, by 2.11(b) and 3.9. Thus: 
$$ 
(Q_{i-1}T)^{f_i}=(Q_{i-1})^{x_i,f_i}T=Q_iT\leq S_{f_{i+1}}, 
$$
and $v\in\bold D$ via $Q_0T$. Now $\Pi(u)\in\ca N\Pi(v)$ by 3.4. 
\qed
\enddemo 

Let $\equiv$ be the equivalence relation on $\ca L$ defined by the partition 
in 3.14(d). In view of  3.14 we may refer 
to the $\equiv$-classes as the {\it maximal cosets} of $\ca N$ in $\ca L$. 

\proclaim {Lemma 3.15} Let $\ca H$ be a partial subgroup of the locality 
$\ca L$, containing the partial normal subgroup $\ca N\norm\ca L$. Then 
$\ca H$ is the disjoint union of the maximal cosets of $\ca N$ contained in $\ca H$. 
\endproclaim 

\demo {Proof} Let $f\in\ca H$. Apply the Frattini lemma (3.11) to obtain 
$f=xh$ for some $x\in\ca N$ and some $h\in N_{\ca L}(T)$ such that $h$ is 
$\up$-maximal with respect to $\ca N$. Then $h=x\i f$ by 1.3(d), and thus 
$h\in\ca H$ as $\ca N\leq\ca H$. Then also $\ca Nh\sub\ca H$, where 
$\ca Nh$ is a maximal coset of $\ca N$ by 3.14(b). 
\qed 
\enddemo

The set $\ca L/\ca N$ of maximal cosets of $\ca N$ may also be denoted 
$\bar{\ca L}$. Let $\r:\ca L\to\bar{\ca L}$ be the mapping which sends 
$g\in\ca L$ to the unique maximal coset of $\ca N$ containing $g$. Set 
$\bold W:=\bold W(\ca L)$ and $\bar{\bold W}=\bold W(\bar{\ca L})$, and let 
$\r^*:\bold W\to\bar{\bold W}$ be the induced mapping of free monoids. 
For any subset or element $X$ of $\bold W$, write $\bar X$ for the image of 
$X$ under $\r^*$, and similarly if $Y$ is a subset or element of $\ca L$ 
write $\bar Y$ for the image of $Y$ under $\r$. In particular, 
$\bar{\bold D}$ denotes the image of $\bold D$ under $\r^*$. Set 
$\bar{\D}=\{\bar P\mid P\in\D\}$. 

For $w\in\bold W$, we shall say that $w$ is $\up$-maximal if every entry 
of $w$ is $\up$-maximal.

\proclaim {Lemma 3.16} There is a unique mapping $\bar\Pi:\bar{\bold D}\to\bar{\ca L}$, a unique 
involutory bijection $\bar g\maps{\bar g}\i$ on $\bar{\ca L}$, and a unique element $\bar{\1}$ 
of $\bar{\ca L}$ such that $\bar{\ca L}$, with these structures, is a partial group, and such that 
$\r$ is a homomorphism of partial groups. Moreover, the homomorphism of free monoids 
$\r^*:\bold W(\ca L)\to\bold W(\bar{\ca L})$ maps $\bold D(\ca L)$ onto $\bold D(\bar{\ca L})$. 
\endproclaim 

\demo {Proof} Let $u=(g_1,\cdots,g_n)$ and $v=(h_1,\cdots,h_n)$ be members of $\bold D$ such that 
$\bar u=\bar v$. By 3.14(d) there exists, for each $i$, an $\up$-maximal $f_i\in\ca L$ with 
$g_i,h_i\in\ca Nf_i$. Set $w=(f_1,\cdots,f_n)$. Then $w\in\bold D$ by 
3.14(e), and then 3.3(a) shows that $\Pi(u)$ and $\Pi(v)$ are elements of 
$\ca N\Pi(w)$. Thus $\bar{\Pi(u)}=\bar{\Pi(w)}=\bar{\Pi(v)}$, and there 
is a well-defined mapping $\bar\Pi:\bar{\bold D}\to\bar{\ca L}$ given by 
$$ 
\bar\Pi(w)=\bar{\Pi(w)}.\tag* 
$$ 
For any subset $X$ of $\ca L$ write $X\i$ for the set of inverses of elements 
of $X$. For any $f\in\ca L$ we then have $(\ca Nf)\i=f\i\ca N\i$ by 
1.1(4). Here $\ca N\i=\ca N$ as $\ca N$ is a partial group, and then 
$(\ca Nf)\i=\ca Nf\i$ by 3.14(a). The inversion map $\ca Nf\maps\ca Nf\i$ is 
then well-defined, and is an involutory bijection on $\bar{\ca L}$. Set $\bar{\1}=\ca N$. 

We now check that the axioms in 1.1, for a partial group, are satisfied by 
the above structures. Since $\bar{\bold D}$ is the image of $\bold D$ under 
$\r^*$, we get $\bar{\ca L}\sub\bar{\bold D}$. Now let 
$\bar w=\bar u\circ\bar v\in\bar{\bold D}$, let $u,v$ be $\up$-maximal 
pre-images in $\bold W$ of $\bar u$, and $\bar v$, and set $w=u\circ v$. Then 
$w$ is $\up$-maximal, and so $w\in\bold D$ by 3.14(e). Then $u$ and $v$ are in $\bold D$, 
and so $\bar u$ and $\bar v$ are in $\bar{\bold D}$. Thus 
$\bar{\bold D}$ satisfies 1.1(1). Clearly, (*) implies that $\bar\Pi$ 
restricts to the identity on $\bar{\ca L}$, so $\bar\Pi$ satisfies 1.1(2). 

Next, let $\bar u\circ\bar v\circ\bar w\in\bar{\bold D}$, and choose corresponding 
$\up$-maximal pre-images $u,v,w$. Set $g=\Pi(v)$. Then 
$\bar g=\bar\Pi(\bar v)$ by (*). By 1.1(3) we have both 
$u\circ v\circ w$ and $u\circ(g)\circ w$ in $\bold D$, and these two words 
have the same image under $\Pi$. Applying $\r^*$ we obtain words in 
$\bar{\bold D}$ having the same image under $\bar\Pi$, and thus $\bar\Pi$ 
satisfies 1.1(3). By definition, $\bar\Pi(\nset)=\bar{\1}$, and then the 
condition 1.1(4) is readily verified. Thus, $\bar{\ca L}$ is a partial 
group. 

By definition, $\bar{\bold D}$ is the image of $\bold D$ under $\r^*$. So, 
in order to check that $\r$ is a homomorphism of partial groups it suffices 
to show that if $w\in\bold D$ then $\bar\Pi(w\r^*)=\Pi(w)\r$. But this is simply the statement (*). 
Moreover, it is this observation which establishes that the given partial 
group structure on $\bar{\ca L}$ is the unique one for which $\r$ is a 
homomorphism of partial groups. We have $f\in Ker(\r)$ if and only if 
$f\g=\bar{\1}=\ca N$. Since $\ca Nf\sub\ca N$ implies $f\in\ca N$, 
and since $\ca N$ is the maximal coset of $\ca L$ containing $\1$, we 
obtain $Ker(\r)=\ca N$. 
\qed 
\enddemo

\vskip .2in 
\noindent 
{\bf Section 4: Quotient localities} 
\vskip .1in 

We continue the setup in which $(\ca L,\D,S)$ is a fixed locality and $\ca N\norm\ca L$ is a partial normal 
subgroup. We have seen in 3.16 that the set $\ca L/\ca N$ of maximal cosets of $\ca N$ inherits 
from $\ca L$ a partial group structure via the projection map $\r:\ca L\to\ca L/\ca N$. The aim now is to 
go further, and to show that $\ca L/\ca N$ is a locality. The argument for this involves some 
subleties: the main problem lies in showing that $\bold D(\ca L/\ca N)$ contains $\bold D_{\bar{\D}}$ 
(see 2.1), where $\bar{\D}$ is the the set of all $P\r$ with $P\in\D$. The following three lemmas are 
intended as steps toward addressing this point.

\proclaim {Lemma 4.1} Let $(\ca L,\D,S)$ be a locality and let $\ca N\norm\ca L$. Then 
$(\ca NS,\D,S)$ is a locality. 
\endproclaim 

\demo {Proof} By 2.9 $\ca NS$ is a partial subgroup of $\ca L$. One observes that $\bold D(\ca NS)$ is 
the subset $\bold D_\D$ of $\bold W(\ca NS)$, as defined in (2.1), and this suffices to show that 
$(\ca NS,\D)$ is objective. As $S$ is a maximal $p$-subgroup of $\ca NS$ there is nothing more that  
needs to be shown. 
\qed 
\enddemo

\proclaim {Lemma 4.2} Let $P$ be a subgroup of $S$, let $\G$ be a non-empty set of $S$-conjugates of $P$, 
and set $X=\bigcup\G$. Assume that $P^x\in\G$ for all $x\in X$. Then either 
$\G=\{P\}$ or $N_S(P)\cap X\nsub P$. 
\endproclaim 

\demo {Proof} Let $\S$ be the set of overgroups $Q$ of $P$ in $S$ such $X\cap Q=P$. Thus $P\in\S$. Regard 
$\S$ as a poset via inclusion, and let $Q$ be maximal in $\S$. If $Q=S$ then $X=P$ and $\G=\{P\}$. On the 
other hand, suppose that $Q\neq S$. Then $Q$ is a proper subgroup of $N_S(Q)$, and the maximality of $Q$ 
implies that there exists $x\in N_S(Q)\cap X$ with $x\notin P$. Since $P^x\leq Q$, and since $P^x\sub X$ 
by hypothesis, we conclude that $P^x=P$. Thus $N_S(P)\cap X\nsub P$. 
\qed 
\enddemo

\proclaim {Theorem 4.3} Let $(\ca L,\D,S)$ be a locality, let $\bar{\ca L}$ be a partial group, and 
let $\b:\ca L\to\bar{\ca L}$ be a homomorphism of partial groups such that the induced map 
$\b^*:\bold W(\ca L)\to\bold W(\bar{\ca L})$ sends $\bold D(\ca L)$ onto $\bold D(\bar{\ca L})$. Set 
$\ca N=Ker(\b)$ and $T=S\cap\ca N$. Further, set $\bold D=\bold D(\ca L)$, 
$\bar{\bold D}=\bold D(\bar{\ca L})$, $\bar S=S\b$, and $\bar\D=\{P\b\mid P\in\D\}$. Then 
$(\bar{\ca L},\bar\D,\bar S)$ is a locality. Moreover: 
\roster 

\item "{(a)}" The fibers of $\b$ are the maximal cosets of $\ca N$. 

\item "{(b)}" For each $\bar w\in\bold W(\bar{\ca L})$ there exists $w\in\bold W(\ca L)$ such that 
$\bar w=w\b^*$ and such that each entry of $w$ is $\up$-maximal relative to the partial normal subgroup 
$\ca N$ of $\ca L$. For any such $w$ we then have $w\in\bold D$ if and only if $\bar w\in\bar{\bold D}$. 

\item "{(c)}" Let $P,Q\in\D$ with $T\leq P\cap Q$. Then $\b$ restricts to a surjection 
$N_{\ca L}(P,Q)\to N_{\bar{\ca L}}(P\b,Q\b)$, and to a surjective homomorphism if $P=Q$.  

\item "{(d)}" $\b$ is an isomorphism if and only if $\ca N=1$.  

\item "{(e)}" We have $(S_g)\b=\bar S_{g\b}$ for each $g\in\ca L$ such that $g$ is $\up$-maximal with 
respect to $\ca N$. 

\endroster 
\endproclaim 

\demo {Proof} The hypothesis that $\bold D\b^*=\bar{\bold D}$ implies that $\b^*$ maps the set of words of 
length 1 in $\ca L$ onto the set of words of length 1 in $\bar{\ca L}$. Thus $\b$ is surjective.  

Let $M$ be a subgroup of $\ca L$. The restriction of $\b$ to $M$ is then a homomorphism of partial groups, 
and hence a homomorphism $M\to M\b$ of groups by 1.13. In particular, $\bar S$ is a $p$-group, and $\bar\D$ 
is a set of subgroups of $\bar S$. 

We have $\ca N\norm\ca L$ by 1.14. Let $\L$ be the set of elements $g\in\ca L$ such that $g$ is 
$\up$-maximal relative to $\ca N$. For any $g\in\L$, $\b$ is constant on the maximal coset $\ca Ng$
(see 3.14) of $\ca N$, so $\b$ restricts to a surjection of $\L$ onto $\bar{\ca L}$. This shows that $\b^*$ 
restricts to a surjection of $\bold W(\L)$ onto $\bold W(\bar{\ca L})$. If $\bar w\in\bar{\bold D}$ then 
there exists $w\in\bold D$ with $w\b^*=\bar w$, and then 3.14(e) shows that such a $w$ may be chosen to be in 
$\bold W(\L)$. Set $\bold D(\L)=\bold D\cap\bold W(\L)$. Thus: 
\roster 

\item "{(1)}" $\b^*$ maps $\bold D(\L)$ onto $\bar{\bold D}$.

\endroster 
Let $\bar w\in\bar{\bold D}$, let $w\in\bold D(\L)$ with $w\b^*=\bar w$, let $\bar a,\bar b\in\bar S$, 
and let $a,b\in S$ with $a\b=\bar a$ and $b\b=\bar b$. Then $a,b\in\L$ by 3.7(a), and 
$(a)\circ w\circ(b)\in\bold D(\L)$ by 2.9. Then $(\bar a)\circ \bar w\circ(\bar b)\in\bar{\bold D}$. 
This shows: 
\roster 

\item "{(2)}" $\bar{\bold D}$ is a $\bar S$-biset (as in 2.9). 

\endroster 
Let $\bar a\in\bar S$, let $a\in S$ be a preimage of $\bar a$, and let $h\in\ca L$ be any preimage of 
$\bar a$. Then $(a,h\i)\in\bold D$, and $ah\i\in\ca N$, so: 
\roster 

\item "{(3)}" The $\b$-preimage of an element $\bar a\in\bar S$ is a maximal coset $\ca Na$, where $a\in S$. 

\endroster 
Fix $\bar g\in\bar{\ca L}$, let $g\in\L$ with $g\b=\bar g$, set $P=S_g$, and set 
$$ 
\bar S_{\bar g}=\{\bar x\in\bar S\mid \bar x^{\bar g}\in \bar S\}.   
$$ 
Let $\bar a\in\bar S_{\bar g}$ and set $\bar b=\bar a^{\bar g}$. As in the proof of 2.6 we may then show 
that $\bar P^{\bar a}\sub\bar S_{\bar g}$. Namely, from $(\bar g\i,\bar a,\bar g)\in\bar{\bold D}$ 
and $\bar{\Pi}(\bar g\i,\bar a,\bar g)=\bar b$ we obtain (from two applications of (2)):
$$ 
(\bar P^{\bar a})^{\bar g}=\bar P^{\bar g\bar b}\leq\bar S. \tag*
$$ 
Thus, the set $\G$ of all $\bar S_{\bar g}$-conjugates of $\bar P$ is a set of subgroups of the set 
$\bar S_{\bar g}$. Setting $X=\bigcup\G$ we thus have the setup of lemma 4.2, with $\bar S_{\bar g}$ in 
the role of $R$. 

Assume now that $\bar P\neq\bar S_{\bar g}$. Then 4.2 yields an element $\bar x\in X-\bar P$ such that 
$\bar x$ normalizes $\bar P$. Let $Q$ be the $\b$-preimage of $\bar P\<\bar x\>$ in $S$. Then $Q^g$ is 
defined (and is a subgroup of $N_{\ca L}(P^g)$) by 2.3(b). As $\bar Q^{\bar g}\leq\bar S$ we obtain 
$Q^g\leq\ca NS$. As $\ca NS$ is a locality by 4.1, it follows from 2.11(b) that there exists $f\in\ca N$ 
with $(Q^g)^f\leq S$. Here $(g,f)\in\bold D$ via $P$, so $Q\leq S_{gf}$, and this contradicts the 
$\up$-maximality of $g$. We conclude: 
\roster 

\item "{(4)}" $(S_g)\b=\bar S_{g\b}$ for each $g\in\L$. 

\endroster 
Thus (e) holds. 

For $\bar w\in\bold W(\bar{\ca L})$ define $\bar S_{\bar w}$ to be the set of all $\bar x\in\bar S$ such that 
$\bar x$ is conjugated successively into $\bar S$ by the entries of $\bar w$. An immediate consequence of 
(4) is then: 
\roster 

\item "{(5)}" $(S_w)\b=\bar S_{w\b^*}$ for all $w\in\bold W(\L)$. 

\endroster 
Notice that (5) implies  point (b). 

We may now verify that $(\bar{\ca L},\bar\D)$ is objective. Thus, let $\bar w\in\bold W(\bar{\ca L})$ with 
$\bar S_{\bar w}\in\bar\D$. Let $w\in\bold W(\L)$ with $w\b^*=\bar w$. Then (5) yields $S_w\in\D$, so 
$w\in\bold D$, and hence $\bar w\in\bar{\bold D}$. Thus $(\bar{\ca L},\bar\D)$ satisfies the condition (O1) 
in definition 2.1 of objectivity. Now let $\bar P\in\bar{\D}$, let $\bar f\in N_{\bar{\ca L}}(\bar P,\bar S)$, 
and set $\bar Q=(\bar P)^{\bar f}$. Then (4) yields $\bar Q=Q\b$ for some $Q\in\D$, and thus 
$\bar Q\in\bar\D$. Any overgroup of $\bar Q$ in $\bar S$ is the image of an overgroup 
of $Q$ in $S$ as $\b$ maps $S$ onto $\bar S$, so $(\bar{\ca L},\bar\D)$ satisfies (O2) in 2.1. Thus 
$(\bar{\ca L},\bar\D)$ is objective. 

Let $P,Q\in\D$ with $T\leq P\cap Q$, set $\bar P=P\b$ and $\bar Q=Q\b$, and let $\bar g\in\bar{\ca L}$ 
such that $\bar P^{\bar g}$ is defined and is a subset of $\bar Q$. Let $g$ be a preimage of $\bar g$ 
in $\L$. Then $P\leq S_g$ by (4). As $\b$ is a homomorphism, $(P^g)\b$ is a subgroup of $\bar Q$. 
Then $P^g\leq Q$ since $\b$ restricts to an epimorphism $S\to\bar S$ with kernel $T$. As $\b$ maps 
subgroups of $\ca L$ homomorphically to subgroups of $\bar{\ca L}$ (by 1.13) we obtain (c). 
In the special case that $P=Q=S$ we obtain in this way an epimorphism from 
$N_{\ca L}(S)$ to $N_{\bar{\ca L}}(\bar S)$. As $S$ is a Sylow subgroup of $N_{\ca L}(S)$ it follows 
that $\bar S$ is a maximal $p$-subgroup of $\bar{\ca L}$, and so $(\bar{\ca L},\bar\D,\bar S)$ is a 
locality. 

Let $g,h\in\L$ with $\g\b=h\b$. Then $S_g=S_h$ by (4), so $(g,h\i)\in\bold D$, and $(gh\i)\b=\1$. Thus 
$g\in\ca Nh$, and then $\ca Ng=\ca Nh$ by 3.14(b). This yields (a), and it remains only to prove (d). 

If $\b$ is an isomorphism then $\b$ is injective, and $\ca N=\1$. On the other hand, suppose that 
$\ca N=\1$. Then (a) shows that $\b$ is injective, and so $\b$ is a bijection. That
$\b\i$ is then a homomorphism of partial groups is given by (5). Thus (d) holds, and the proof is 
complete.  
\qed 
\enddemo

\definition {Definition 4.4} Let $\ca L$ and $\ca L'$ be partial groups and let 
$\b:\ca L\to\ca L'$ be a homomorphism. Then $\b$ is a {\it projection} if $\bold D\b^*=\bold D'$.  
\enddefinition 

\proclaim {Corollary 4.5} Let $(\ca L,\D,S)$ be a locality, let $\ca N\norm\ca L$ be a partial normal 
subgroup and let $\r:\ca L\to\ca L/\ca N$ be the mapping which sends $g\in\ca L$ to the unique maximal 
coset of $\ca N$ containing $g$. Set $\bar{\ca L}=\ca L/\ca N$, set $\bar S=S\r$, and let $\bar\D$ be 
the set of images under $\r$ of the members of $\D$. Regard $\bar{\ca L}$ as a partial group in the 
unique way (given by 3.16) which makes $\r$ into a homomorphism of partial groups. Then 
$(\bar{\ca L},\bar\D,\bar S)$ is a locality, and $\r$ is a projection. 
\endproclaim 

\demo {Proof} Immediate from 3.16 and 4.3. 
\qed 
\enddemo 

\proclaim {Theorem 4.6 (``First Isomorphism Theorem")} Let $(\ca L,\D,S)$ and $(\ca L',\D',S')$ be 
localities, let $\b:\ca L\to\ca L'$ be a projection, and let $\ca N\norm\ca L$ be a partial normal 
subgroup of $\ca L$ contained in $Ker(\b)$. Let $\r:\ca L\to\ca L/\ca N$ be the 
projection given by 4.5. Then there exists a unique homomorphism 
$$ 
\g:\ca L/\ca N\to\ca L' 
$$ 
such that $\r\circ\g=\b$, and $\g$ is a projection. Moreover, $\g$ is an isomorphism if and only if 
$\ca N=Ker(\b)$. 
\endproclaim  

\demo {Proof} Set $\ca M=Ker(\b)$, and let $h\in\ca L$ be $\up$-maximal relative to $\ca M$. Then 
$\ca Mh$ is a maximal coset of $\ca M$ in $\ca L$ by 3.14(b), and $\ca Mh=\ca Mh\ca M$ by 3.14(a). Let 
$g\in\ca Mh$. As $\ca N\leq\ca M$ by hypothesis, the splitting lemma (3.12, as applied to $\ca M$ and $h$) 
yields 
$$ 
\ca Ng\ca N\sub\ca N(\ca Mh)\ca N\sub\ca M(\ca Mh)\ca M=\ca Mh\ca M. 
$$  
The definition 3.6 of the relation $\up$ on $\ca L\circ\D$ shows that $\ca Ng\ca N$ contains an element 
$f$ which is $\up$-maximal with respect to $\ca N$. Then $f\in\ca Mh\ca M$, so $f\in\ca Mh$, and another 
application of the splitting lemma yields $\ca Nf\sub\ca Mh$. Here $\ca Nf$ is the maximal coset of $\ca N$ 
containing $f$. We have thus shown: 
\roster 

\item "{(*)}" The partition $\ca L/\ca N$ of $\ca L$ is a refinement of the partition $\ca L/\ca M$. 

\endroster 
By 4.3(a) $\b$ induces a bijection $\ca L/\ca M\to\ca L'$.  Set $\bar{\ca L}=\ca L/\ca N$. Then (*) 
implies that there is a mapping $\g:\bar{\ca L}\to\ca L'$ which sends the maximal coset $\ca Nf\ca N$ 
to $f\b$. Clearly, $\g$ is the unique mapping $\bar{\ca L}\to\ca L'$ such that $\r\circ\g=\b$. 

Let $\bar w\in\bold D(\bar{\ca L})$. Then 4.3(b) yields a word $w\in\bold D$ such that $w\r^*=\bar w$ and 
such that the entries of $w$ are $\up$-maximal relative to $\ca N$. We have $\bar w\g^*=w\b^*$, so 
$\g^*$ maps $\bold D(\bar{\ca L})$ into $\bold D(\ca L')$. Let $\Pi'$ and $\bar\Pi$ be the products in 
$\ca L'$ and $\bar{\ca L}$, respectively. As $\b$ and $\r$ are homomorphisms we get 
$$ 
\Pi'(\bar w\g^*)=\Pi'(w\b^*)=(\Pi(w))\b=(\bar{\Pi(w)})\g=(\bar\Pi(\bar w))\g,  
$$ 
and thus $\g$ is a homomorphism. As $\b=\r\circ\g$ is a projection, and $\r$ is a projection, one 
verifies that $\g^*$ maps $\bold D(\bar{\ca L})$ onto $\bold D(\ca L')$ and that $\g$ maps 
$\bar\D$ onto $\D'$. Thus $\g$ is a projection. 

We have $\ca M=\ca N$ if and only if $Ker(\g)=\1$. Then 4.3(d) shows that $\g$ is an isomorphism if 
and only if $\ca M=\ca N$; completing the proof. 
\qed 
\enddemo

\proclaim {Proposition 4.7 (Partial Subgroup Correspondence)} Let $(\ca L,\D,S)$ and 
$(\bar{\ca L},\bar\D,\bar S)$ be localities, and let $\b:\ca L\to\bar{\ca L}$ be a projection. Set 
$\ca N=Ker(\b)$ and set $T=S\cap\ca N$. Then $\b$ induces a bijection $\s$ from the set $\frak H$ of 
partial subgroups $\ca H$ of $\ca L$ containing $\ca N$ to the set $\bar{\frak H}$ of partial subgroups 
$\bar{\ca H}$ of $\bar{\ca L}$. Moreover, for any $\ca H\in\frak H$, we have $\ca H\b\norm\ca L'$ if 
and only if $\ca H\norm\ca L$. 
\endproclaim 

\demo {Proof} Any partial subgroup of $\ca L$ containing $\ca N$ is a union of maximal cosets of $\ca N$ by 
3.15. Then 4.3(a) enables the same argument that one has for groups, for proving that $\r$ induces a 
bijection $\frak H\to\bar{\frak H}$. Since each maximal coset of $\ca N$ contains an element which is 
$\up$-maximal with respect to $\ca N$, one may apply 4.3(b) in order to show that a partial subgroup 
$\ca H\in\frak H$ is normal in $\ca L$ if and only if its image is normal in $\bar{\ca L}$. The 
reader should have no difficulty with the details of the argument. 
\qed 
\enddemo

\definition {4.8 Remark}
A comprehensive ``second isomorphism theorem" appears to be out of reach for two reasons. First, given a 
partial subgroup $\ca H\leq\ca L$ and a partial normal subgroup $\ca N\norm\ca L$, there appears to be no 
reason for the image of $\ca H$ under the projection $\r:\ca L\to\ca L/\ca N$ to be a partial subgroup of 
$\ca L/\ca N$, other than in special cases. Second, there seems to be no way, in general, to define the 
quotient of $\ca H$ over the partial normal subgroup $\ca H\cap\ca N$ of $\ca H$. On the other hand, a
``third isomorphism theorem" may easily be deduced from 4.3  
and from the trivial observation that a composition of projections is again a projection. 
\enddefinition

\proclaim {Lemma 4.9} Let $\ca N\norm\ca L$ and let $\r:\ca L\to\ca L/\ca N$ be the canonical projection. 
Further, let $\ca H$ be a partial subgroup of $\ca L$ containing $\ca N$ and let $X$ be an 
arbitrary subset of $\ca L$. Then $(X\cap\ca H)\r=X\r\cap\ca H\r$. 
\endproclaim 

\demo {Proof} By 3.15, $\ca H$ is a union of maximal cosets of $\ca N$, and then $\ca H\r$ is the set of 
those maximal cosets. On the other hand $X\r$ is the set of all maximal cosets $\ca Ng$ of $\ca N$ such 
that $X\cap\ca Ng\neq\nset$. Thus $X\r\cap\ca H\r\sub(X\cap\ca H)\r$. The reverse inclusion is obvious. 
\qed 
\enddemo 

\proclaim {Corollary 4.10} Let $\ca N\norm\ca L$, and let $\ca M$ be a 
partial normal subgroup of $\ca L$ containing $\ca N$. Let $\r:\ca L\to\ca L/\ca N$ be the canonical 
projection. Then $(S\cap\ca M)\r$ is a maximal $p$-subgroup of $\ca M\r$. 
\endproclaim 

\demo {Proof} Write $(\bar{\ca L},\bar\D,\bar S)$ for the quotient locality given by 4.5, and set 
$\bar{\ca M}=\ca M\r$. Applying 4.9 with $S$ in the role of $X$, we obtain 
$(S\cap\ca M)\r=\bar S\cap\bar{\ca M}$. Since $\bar{\ca M}\norm\bar{\ca L}$, it follows from 2.11(c) that 
$\bar S\cap\bar{\ca M}$ is maximal in the poset of $p$-subgroups of $\bar{\ca M}$, completing the proof. 
\qed 
\enddemo

\proclaim {Proposition 4.11} Let $\ca N\norm\ca L$, set $T=S\cap\ca N$, and set $\ca L_T=N_{\ca L}(T)$. 
Set $\bar{\ca L}=\ca L/\ca N$ and let $\r:\ca L\to\bar{\ca L}$ be the canonical projection. 
Then the partial subgroup $\ca L_T$ of $\ca L$ is a locality $(\ca L_T,\D,S)$, and the  
restriction of $\r$ to $\ca L_T$ is a projection $\ca L_T\to\bar{\ca L}$. 
\endproclaim 

\demo {Proof} That $\ca L_T$ is a partial subgroup of $\ca L$ having the structure of a locality 
$(\ca L_T,\D,S)$ is given by 2.13. Let $\r_T$ be the restriction of $\r$ to $\ca L_T$. Then $\r_T$ is 
a homomorphism of partial groups, and 3.14(e) shows that $\r_T$ maps $\bold D(\ca L_T)$ onto 
$\bold D(\bar{\ca L})$. That is, $\r_T$ is a projection. 
\qed 
\enddemo

We end this section with an application. For $G$ a finite group, $O_{p'}(G)$ denotes the largest normal 
subgroup of $G$ having order prime to $p$, and $G$ is of {\it characteristic $p$} if 
$C_G(O_p(G))\leq O_p(G)$. We will assume that the reader is familiar with the definition of a 
{\it fusion system} over a finite $p$-group. For $(\ca L,\D,S)$ a loclality one has the fusion 
system $\ca F_S(\ca L)$ on $S$, whose isomorphisms are the conjugation maps $c_w$ (for $w\in\bold W(\ca L)$) 
that were introduced following 2.16.

\proclaim {Proposition 4.12} Let $(\ca L,\D,S)$ be a locality. For each $P\in\D$ set 
$\Theta(P)=O_{p'}(N_{\ca L}(P))$, and set $\Theta=\bigcup\{\Theta(P)\}_{P\in\D}$. Assume: 
\roster 

\item "{(*)}" $P\in\D\implies N_{\ca L}(P)/\Theta(P)$ is of characteristic $p$. 

\endroster 
Then $\Theta\norm\ca L$, $S\cap\Theta=1$, and the canonical projection $\r:\ca L\to\ca L/\Theta$ resricts 
to an isomorphism $S\to S\r$. Moreover, upon identifying $S$ with $S\r$: 
\roster 

\item "{(a)}" $(\ca L/\Theta,\D,S)$ is a locality.  

\item "{(b)}" $\ca F_S(\ca L/\Theta)=\ca F_S(\ca L)$. 

\item "{(c)}" For each $P\in\D$, the restriction 
$$ 
\r_P:N_{\ca L}(P)\to N_{\ca L/\Theta}(P)
$$ 
of $\r$ induces an isomorphism 
$$ 
N_{\ca L/\Theta}(P)\cong N_{\ca L}(P)/\Theta(P), 
$$ 
and $N_{\ca L/\Theta}(P)$ is of characteristic $p$. 

\endroster 
\endproclaim 

\demo {Proof} Let $x\in\Theta$. Then there exists $Q\in\D$ with $x\in\Theta(Q)$. Choose such a $Q$ so 
that $|Q|$ is as large as possible, and set $R=N_{S_x}(Q)$. Then $[R,x]\leq RR^x\leq S$. But also 
$$ 
[R,x]\leq [N_{\ca L}(Q),\Theta(Q)]\leq\Theta(Q), 
$$ 
and so $[R,x]=1$. Then $x\in\Theta(R)$ by (*) (with $R$ in the role of $Q$), and the maximality of 
$|Q|$ yields $Q=R$. Thus $Q=S_x$, and we have thus shown that $x\in\Theta(S_x)$. Now let $P\in\D$ with 
$P\leq S_x$. Then $[P,x]\leq [S_x,x]=1$, and (*) yields $x\in\Theta(P)$. Thus: 
\roster 

\item "{(**)}" Let $x\in\Theta$, and let $P\in\D$ with $P\leq S_x$. Then $x\in\Theta(P)$. 

\endroster 

Clearly, $\1\in\Theta$, and $\Theta$ is closed under inversion. Let 
$$ 
w=(x_1,\cdots,x_n)\in\bold W(\Theta)\cap\bold D, 
$$ 
and set $P=S_w$. By (**), and by induction on $n$, we 
obtain $x_i\in\Theta(P)$ for all $i$, and hence $\Pi(w)\in\Theta(P)$. Thus $\Theta$ is a partial subgroup 
of $\ca L$. Now let $x\in\Theta$ and let $g\in\ca L$ be given such that $(g\i,x,g)\in\bold D$ via 
some $Q\in\D$. Then $Q^{g\i}\leq S_x$, so (**) yields $x\in\Theta(Q^{g\i})$, and then 
$x^g\in\Theta(Q)$ by 2.3(b). This completes the proof that $\Theta\norm\ca L$. 

Set $\bar{\ca L}=\ca L/\Theta$ and adopt the usual ``bar"-convention for images of elements, 
subgroups, and collections of subgroups under the quotient map $\r:\ca L\to\bar{\ca L}$. 
Since $\Theta$ is a set of $p'$-elements of $\ca L$ we have $S\cap\Theta=\1$, and we may 
therefore identify $S$ with $\bar S$, and $\D$ with $\bar\D$. Point (a) is then given by 4.5. 

For each $P\in\D$ let $\r_P$ be the restriction of $\r$ to $N_{\ca L}(P)$. 
Then $\r_P$ is an epimorphism $N_{\ca L}(P)\to N_{\bar{\ca L}}(P)$ by 4.3(c), with kernel 
$\Theta(P)$. This yields point (c). 

By 4.3(c) the conjugation maps $c_g:P\to Q$ in $\ca F$, with $P,Q\in\D$ and with $g\in\ca L$, are the same 
as the conjugation maps $c_{\bar g}:P\to Q$ with $\bar g\in\ca L/\Theta$. Since $\ca F_S(\ca L)$ is 
$\D$-generated (by definition 2.12), we obtain $\ca F_S(\ca L)=\ca F_S(\ca L/\Theta)$. That is, (b) 
holds, and the proof is complete. 
\qed 
\enddemo

\vskip .1in 
\noindent 
{\bf Section 5: Products of partial normal subgroups} 
\vskip .1in 

There are two main results in this section. The first (Theorem 5.1) concerns products of partial 
normal subgroups in a locality. The second (Proposition 5.6) is an application of essentially all of 
the results preceding it, and will play a vital role in the successor to this work.

\proclaim {Theorem 5.1} Let $(\ca L,\D,S)$ be a locality, and let $\ca M\norm\ca L$ and $\ca N\norm\ca L$ 
be partial normal subgroups. Set $U=S\cap\ca M$ and $V=S\cap\ca N$, and assume: 
\roster 

\item "{(*)}" $\ca M$ normalizes $V$, and $\ca N$ normalizes $U$. 

\endroster 
Then $\ca M\ca N=\ca N\ca M\norm\ca L$, and $S\cap\ca M\ca N=UV$. 
\endproclaim 

The proof will require the following version of the Splitting Lemma (3.12). 

\proclaim {Lemma 5.2} Assume the hypothesis of 5.1, and let $g\in\ca M\ca N$. Then there exists 
$(x,y)\in\bold D$ with $x\in\ca M$, $y\in\ca N$, $g=xy$, and $S_g=S_{(x,y)}$. 
\endproclaim 

\demo {Proof} Consider the set of all triples $(g,x,y)\in\ca M\ca N\times\ca M\times\ca N$ such that 
$g$ is a counter-example to the lemma, and $g=xy$. Among all such triples, let $(g,x,y)$ be chosen 
so that $|S_{(x,y)}|$ is as large as possible. Set $Q=S_{(x,y)}$ and set $P=N_{S_g}(Q)$. It suffices 
to show that $P=Q$ in order to obtain the lemma. 

By 3.2 we have $(y,y\i,x,y)\in\bold D$ and $g=yx^y$, with $S_{(x,y)}=S_{(y,x^y)}$. Suppose that 
$P\leq S_y$. Then $P^y\leq S$, and since $P^g\leq S$ we conclude that 
$P\leq S_{(y,x^y)}$, and hence $P=Q$, as desired. Thus we may assume: 
\roster 

\item "{(1)}" $P\nleq S_y$. 

\endroster 
Let $h$ be $\up$-maximal in the maximal coset of $\ca M$ containing $g$. Then 3.9 yields an element 
$r\in\ca M$ such that $g=rh$, and 3.12 yields $S_g=S_{(r,h)}$. Then $Q\leq S_{(r,h)}$, so 
$(y\i,x\i,r,h)\in\bold D$ via $Q^g$ and $\Pi(y\i,x\i,r,h)=\Pi(g\i,g)=\1$. Thus: 
$$ 
y=x\i rh\quad\text{and}\quad h=r\i xy. \tag* 
$$ 
Since $y,h\in N_{\ca L}(U)$, it follows that $r\i x\in N_{\ca M}(U)$, and 
then that $h=(r\i x)y\in\ca M\ca N$. 

Suppose that $h$ does not provide a counter-example to the lemma. That is, suppose that there exists 
$x'\in\ca M$ and $y'\in\ca N$ such that  $(x',y')\in\bold D$, $x'y'=h$, and $S_{(x',y'})=S_h$. As 
$r\i xy=h=x'y'$ we get $xy=rx'y'$, and $(rx',y')\in\bold D$ with $rx'y'=rh=g$. The idea now is to 
replace $(x,y)$ with $(rx',y')$ and to contradict the assumption that $S_g\neq Q$. In order to achieve this, 
observe first of all that $S_g\leq S_r$ since $S_{(r,h)}=S_{rh}=S_g$. Then observe that 
$(S_g)^r\leq S_h$, and that $S_h=S_{(x',y')}\leq S_{x'}$. Thus $(S_g)^r\leq S_{x'}$, and so 
$S_g\leq S_{rx'}$. As $rx'y'=g$ we conclude that $S_g\leq S_{(rx',y')}$, which yields the desired 
contradiction. We conclude that $h$ is itself a counter-example to the lemma. 

Since $h=r\i xy$ by (*), and since $h$ and $y$ are in $N_{\ca L}(U)$, we have $r\i x\in N_{\ca M}(U)$, 
and then $U\leq S_{(r\i x,y)}$ since $h\in N_{\ca L}(U)$. Note furthermore that 
$Q=S_{(x,y)}\leq S_g=S_{(r,h)}\leq S_r$, and thus $Q^rU\leq S_{(r\i x,y)}$. The maximality of 
$|Q|$ in our initial choice of $(g,x,y)$ then yields $Q^r=Q^rU=S_{(r\i x,y)}$. Thus $U^r\leq Q^r$, 
and conjugation by $r\i$ yields $U\leq Q$. A symmetric argument yields $V\leq Q$. 
Setting $H=N_{\ca L}(Q)$, it now follows from 3.1(b) that $x,y\in H$. Then $Q=O_p(H)$. 

Set $X=H\cap\ca M$ and $Y=H\cap\ca N$. Then $X,Y$, and $UV$ are normal 
subgroups of $H$, and $XY/UV$ is a $p'$-group. Set $\bar H=H/(X\cap Y)UV$. 
Here $P\leq H$ and $[P,g]\leq S$. Since $g\in XY$ we obtain 
$$
[\bar P,\bar g]=[\bar P,\bar x][\bar P,\bar y]\leq \bar X\bar Y,  
$$  
and since $\bar X\bar Y$ is a $p'$-group we get $[\bar P,\bar g]=1$. 
As $\bar X\cap\bar Y=1$ we have  
$C_{\bar X\bar Y}(\bar P)=C_{\bar X}(\bar P)\times C_{\bar Y}(\bar P)$. 
As $\bar g=\bar x\bar y$ it follows that $\bar x$ and $\bar y$ centralize 
$\bar P$. Thus $P^x\leq(X\cap Y)P$ and $P\in Syl_p((X\cap Y)P)$. By Sylow's 
Theorem there exists $z\in X\cap Y$ with $P^x=P^z$. Replacing $(x,y)$ with 
$(xz\i,zy)$ we get $g=(xz\i)(zy)$ and $P\leq S_{(xz\i,zy)}$. This 
contradicts the maximality of $Q$ and yields a final contradiction, proving the lemma. 
\qed 
\enddemo 

\demo {Proof of 5.1} Let $w=(g_1,\cdots,g_n)\in\bold W(\ca M\ca N)\cap\bold D$ via $Q\in\D$. By 4.11 we 
may write $g_i=x_iy_i$ with $x_i\in\ca M$, $y_i\in\ca N$, and with $S_{g_i}=S_{(x_i,y_i)}$. Set 
$w'=(x_1,y_1,\cdots,x_n,y_n)$. Then $w'\in\bold D$ via $Q$ and $\Pi(w)=\Pi(w')$. Since each $y_i$ normalizes 
$U$, it follows from 3.4 that $\Pi(w')=\Pi(w'')$ for some $w''$ such that $w''=u\circ v\in\bold D$, where 
$u\in\bold W(\ca M)$, and where $v\in\bold W(\ca N)$. Thus $\ca M\ca N$ is closed under $\Pi$. 
In order to show that $\ca M\ca N=(\ca M\ca N)\i$ we note that if $(x,y)\in\bold D\cap(\ca M\times\ca N)$ 
then $(y\i,x\i)\in\bold D$ and that $y\i x\i\in\ca M\ca N$ by 3.2. Thus $\ca M\ca N$ is a partial 
subgroup of $\ca L$. Moreover, we have shown that $\ca M\ca N=\ca N\ca M$.

Let $g\in\ca M\ca N$ and let $f\in\ca L$ with $(f\i,g,f)\in\bold D$. As usual we may write $f=hr$ with 
$r\in\ca N$, $h\in N_{\ca L}(V)$, and $S_f=S_{(h,r)}$. Write $g=xy$ as in 5.2. By assumption we have 
$(f\i,g,f)\in\bold D$ via some $P\in\D$. Setting $v=(r\i,h\i,x,y,h,r)$ it follows that 
$v\in\bold D$ via $P$ and that $g^f=\Pi(v)$. Here $(h\i,h,y,h)\in\bold D$ via $S_{(y,h)}$ by 3.2, so 
$v':=(r\i,h\i,x,h,h\i,y,h,r)\in\bold D$ via $P$. Then  
$$
g^f=\Pi(v)=\Pi(v')=(x^h y^h)^r\in(\ca M\ca N)^r. 
$$ 
Since $r\in\ca N$, and $\ca M\ca N$ is a partial group, we conclude that 
$g^f\in\ca M\ca N$. Thus $\ca M\ca N\norm\ca L$. 

Set $M=N_{\ca M}(S)$, $N=N_{\ca N}(S)$, and let $s\in S\cap\ca M\ca N$. Then 5.2 yields $s=fg$ with 
$f\in\ca M$, $g\in\ca N$, and with $S=S_{(f,g)}$. Thus $f\in M$ and $g\in N$, where $M$ and $N$ are 
normal subgroups of the group $N_{\ca L}(S)$. Then $UV$ is a normal Sylow $p$-subgroup of $MN$, and since 
$s=fg\in MN$ we obtain $s\in UV$. Thus $S\cap\ca M\ca N=UV$, and the proof is complete. 
\qed 
\enddemo

\proclaim {Lemma 5.3} Let $(\ca L,\D,S)$ be a locality, let $\ca M$ and $\ca N$ be partial normal 
subgroups of $\ca L$, and set $U=S\cap\ca M$ and $V=S\cap\ca N$. Suppose that $\ca M\cap\ca N\leq S$. 
Then $\ca M\leq N_{\ca L}(V)$ and $\ca N\leq N_{\ca L}(U)$. 
\endproclaim 

\demo {Proof} Let $g\in\ca M$, set $P=S_g$, and let $x\in N_V(P)$. Then 
$(x\i,g\i,x,g)\in\bold D$ via $P^{gx}$, and then 
$x\i g\i xg\in\ca M\cap\ca N$. The hypothesis then yields $x^g\in S$, and 
thus $N_V(P)\leq P$. Since $PV$ is a subgroup of $S$ we conclude that 
$V\leq P$, and then $V^g=V$ by 3.1(a). Thus $\ca M \leq N_{\ca L}(V)$, 
and the lemma follows. 
\qed 
\enddemo

\proclaim {Corollary 5.4} Let $\ca M,\ca N\norm\ca L$ and suppose that $\ca M\cap\ca N\leq S$. 
Then $\ca M\ca N\norm\ca L$, and $S\cap\ca M\ca N=(S\cap\ca M)(S\cap\ca N)$. 
\endproclaim 

\demo {Proof} Immediate from 5.1 and 5.3. 
\qed 
\enddemo

The proof of the following result (which will play an essential role in Part III) uses essentially 
everything that has preceded it.

\proclaim {Proposition 5.5} Let $(\ca L,\D,S)$ be a locality, let $\ca N\norm\ca L$ be a partial 
normal subgroup, set $T=S\cap\ca N$, set $\ca L_T=N_{\ca L}(T)$, and let $\ca K$ be a partial normal 
subgroup of $\ca L_T$. Assume: 
\roster 

\item "{(1)}" $C_{N_{\ca L}(P)}(O_p(N_{\ca L}(P))\leq O_p(N_{\ca L}(P)$ for all $P\in\D$, and 

\item "{(2)}" $\ca K\leq C_{\ca L}(T)$. 

\endroster 
Then $\<\ca K,\ca N\>\norm\ca L$, and $S\cap\<\ca K,\ca N\>=(S\cap\ca K)T$. Moreover, if 
$S=C_S(T)T$ then $\<\ca K,\ca N\>=\ca K\ca N=\ca N\ca K$. 
\endproclaim 

\demo {Proof} Set $\ca N_T=N_{\ca N}(T)$. Then $\ca N_T\norm\ca L_T$ by 1.8. View $\ca L_T$ as a 
locality $(\ca L_T,\D,S)$ is in 2.13. The condition (1) allows us to apply 3.7, and to thereby conclude 
that $\ca N_T$ normalizes $C_S(T)T$. Set $U=S\cap\ca K$. Then $U\leq C_S(T)$ by (2), and then 
$[U,\ca N_T]\leq C_S(T)T\cap\ca K=U$. Thus $\ca N_T$ normalizes $U$. Since $\ca K$ normalizes $T$, we 
conclude from 5.1 that $\ca K\ca N_T\norm\ca L_T$.  

Let $\bar{\ca L}$ be the quotient locality $\ca L/\ca N$ and let $\r:\ca L\to\bar{\ca L}$ be the 
canonical projection. Then the restriction of $\r$ to $\ca L_T$ is a projection $\ca L_T\to\bar{\ca L}$ 
by 4.11, and $(\ca K\ca N_T)\r=\ca K\r$. Then $\ca K\r\norm\bar{\ca L}$ by partial subgroup correspondence 
(4.7). The $\r$-preimage of $\ca K\r$ in $\ca L$ is then a partial normal subgroup of $\ca L$, 
containing the partial subgroup $\<\ca K,\ca N\>$ of $\ca L$ generated by $\ca K$ and $\ca N$. 

By 1.9, $\<\ca K,\ca N\>$ is the union of its subsets $Y_i$, where $Y_0=\ca K\cup\ca N$ and (for 
$k>0$) $Y_k$ is the set of all $\Pi(w)$ with $w\in\bold W(Y_{k-1})\cap\bold D$. Clearly $\r$ maps 
$Y_0$ into $\ca K\r$, and a straight-forward induction on $k$ then shows that $\r$ maps each $Y_k$ 
into $\ca K\r$. Thus $\<\ca K,\ca N\>$ is mapped onto $\ca K\r$, and partial subgroup correspondence 
then implies that $\<\ca K,\ca N\>$ is the preimage of $\ca K\r$. As $\ca K\r\norm\bar{\ca L}$, 
a further application of partial subgroup correspondence yields $\<\ca K,\ca N\>\norm\ca L$. 

Set $U=S\cap\ca K$ and set $V=S\cap\<\ca K,\ca N\>$. Then the restriction of $\r$ to $V$ is a 
homomorphism of groups by 1.16, with kernel $T$. Here $U$ is a maximal $p$-subgroup of $\ca K$ 
by 3.1(c), and both $V\r$ and $U\r$ are maximal $p$-subgroups of $\ca K\r$ by 4.10. Then 
$V\r=U\r$ and $V=UT$.  

Suppose now that $S=C_S(T)T$. Then each element of $\ca K$ is $\up$-maximal with respect to $\ca N$, 
by 3.7. Let $w\in\bold W(\ca K\ca N)\cap\bold D$, write $w=(x_1y_1,\cdots,x_ny_n$ with 
$x_i\in\ca K$ and $y_i\in\ca N$, and set $w'=(x_1,y_1,\cdots,x_n,y_n)$. Then $w'\in\bold D$ by the 
Splitting Lemma (3.12), and $\Pi(w')=\Pi(w')$. Here $\Pi(w')\in\ca K\ca N$ by 3.4, 
so $\<\ca K,\ca N\>=\ca K\ca N$ in this case. One similarly has $\<\ca K,\ca N\>=\ca N\ca K$, completing 
the proof. 
\qed 
\enddemo

\definition {Remark 5.6} In the proof of 5.5 hypothesis (2) serves no other purpose than to 
guarantee that $\ca K\ca N_T$ is a partial normal subgroup of $\ca L_T$. In fact, by [Theorem A in He], 
the product of partial normal subgroups of a locality is always a partial normal subgroup, and so 
(2) is redundant. 
\enddefinition

\vskip .2in 
\noindent 
{\bf Appendix A: Limits and colimits in the category of partial groups} 
\vskip .1in 

This appendix was inspired by some remarks of Edoardo Salati, who identified a serious gap in the author's 
earlier treatment of colimits, and who has himself shown [Sal] that the category of 
partial groups is complete (has all limits) and co-complete (has all colimits). The discussion here 
will establish a somewhat weaker result. 

By a {\it pointed set} we mean a set with a distinguished base-point, and there is then a category 
$Set^*$ of pointed sets with base-point-preserving maps. Let $Part$ be the category of partial groups. 
There is then a forgetful functor $Part\to Set^*$, given by regarding a partial group as a pointed set 
having the identity element as its base-point. 

In order to discuss limits and colimits in $Part$ (and their relation with limits and colimits in $Set^*$) 
we begin by reviewing some definitions. 

\definition {Definition A1} Let $J$ be a small category and let $\ca C$ be a category. By a {\it $J$-shaped 
diagram in $\ca C$} we mean a covariant functor $F:J\to\ca C$. 
\enddefinition 

As always, composition of mappings will be written from left to right. 

\definition {Definition A2} Let $F:J\to\ca C$ be a $J$-shaped diagram in $\ca C$. A {\it cone} to $F$ consists 
of an object $M$ of $\ca C$ together with a family $\phi=(\phi_X:M\to F(X))_{X\in Ob(J)}$ of $\ca C$-morphisms, 
such that for each $J$-morphism $f:X\to Y$ we have $\phi_Y=\phi_X\circ F(f)$. The cone 
$(M,\phi)$ is a {\it limit} of $F$ if for every cone $(N,\psi)$ to $F$ there exists a 
unique $\ca C$-morphism $u:N\to M$ such that $\psi_X=u\circ\phi_X$ for all $X\in Ob(J)$. 
\enddefinition

Consider now the case in which $\ca C$ is the category of sets (and mappings of sets), and let $F:J\to\ca C$ 
be a $J$-shaped diagram. If the only $J$-morphisms are identity morphisms then the direct product 
$\wh M$ of the sets $F(X)$ for $X\in Ob(J)$, together with the set $\wh\phi$ of associated projection maps 
is a limit of $F$. More generally,  let $M$ be the subset of $\wh M$ consisting of all $Ob(J)$-tuples 
$(a_X)_{X\in Ob(J)}$ such that, for 
each $J$-morphism $f:X\to Y$, we have $a_Y=(a_X)F(f)$ ($a_Y$ is equal to the image of $a_X$ under $F(f)$). 
Then $M$, together with the set $\phi$ of maps $\phi_X:M\to F(X)$ where $\phi_X$ is the restriction to $M$  
of the projection $\wh\phi_X:\wh M\to F(X)$, is a limit of $F$. If instead $\ca C$ is taken to be the 
category $Set^*$ of pointed sets, then $\wh M$ and $M$ are pointed sets (via the $Ob(J)$-tuple of base-points 
$*_X\in F(X)$) and one observes that $(M,\phi)$ is again a limit of $F$.

\proclaim {Theorem A.3} Let $Part$ be the category of partial groups, let $Set^*$ be the category of pointed 
sets, let $J$ be a small category, and let $F:J\to Part$ be a $J$-shaped diagram.  
Let $F_0:J\to Sets^*$ be the composition of $F$ with the 
forgetful functor $Part\to Sets^*$. Then there exists a limit $(\ca M,\phi)$ of $F$, and the forgetful 
functor $Part\to Set^*$ sends $(\ca M,\phi)$ to a limit of $F_0$. 
\endproclaim 

\demo {Proof} We shall only outline the steps to the proof, leaving most details to the reader. 
Let $\wh{\ca M}$ be the pointed set obtained as the direct product of the partial 
groups $F(X)$ for $X\in Ob(J)$. Let $\wh{\bold D}$ be the direct product of the pointed sets 
$\bold D(F(X))$. Thus the members of $\wh{\bold D}$ are $Ob(J)$-tuples $(w_X)_{X\in Ob(J)}$, 
with $w_X\in\bold D(F(X))$. Let $\Pi_X:\bold D(F(X))\to F(X)$ be the product. There is then a mapping 
$$ 
\wh\Pi:\wh{\bold D}\to\wh{\ca M} 
$$ 
which sends $(w_X)_{X\in Ob(J)}$ to $(\Pi_X(w_X))_{X\in Ob(J)}$. It is now straightforward to check that 
$\wh{\ca M}$ is a partial group via the product $\wh\Pi$ and via the inversion map which sends an 
element $(g_X)$ of $\wh{\ca M}$ to the $Ob(J)$-tuple $(g_X\i)$ of inverses. 

Let $\ca M$ be the subset of $\wh{\ca M}$ consisting of all $Ob(J)$-tuples 
$(g_X)_{X\in Ob(J)}$ such that, for each $J$-morphism $f:X\to Y$, we have $g_Y=(g_X)F(f)$. Let 
$\phi$ be the $Ob(J)$-tuple $(\phi_X)$ of maps $\phi_X:\ca M\to F(X)$ obtained by restriction to $\ca M$ 
of the projection $\wh\phi_X:\wh{\ca M}\to F(X)$. One observes that each $\phi_X$ is a homomorphism 
of partial groups, and that $(\wh{\ca M},\phi)$ is a cone of $F$. 

Now let $(\ca N,\psi)$ be any cone of $F$. For $g\in\ca N$ define $g\mu$ to be the $Ob(J)$-tuple 
$(g\psi_X)$. One checks that each such $g\mu$ is an element of $\ca M$, and then that the mapping 
$\mu:\ca N\to\ca M$ is a homomorphism of partial groups. Finally, one observes that 
$\psi_X=u\circ\phi_X$ for all $X\in Ob(J)$ and that $u$ is necessarily the unique homomorphism 
$\ca N\to\ca N$ having this property. Thus $(\ca M,\phi)$ is a limit of $F$. But also
$(\ca M,\phi)$ is a limit of $F_0$, which completes the proof. 
\qed 
\enddemo

The situation for colimits of partial groups is not as straightforward as that of limits, as we shall see.  

\definition {Definition A.4} Let $F:J\to\ca C$ be a $J$-shaped diagram in $\ca C$. A {\it co-cone} to $F$ 
consists of an object $M$ of $\ca C$ together with a family $\phi=(\phi_X:F(X)\to M)_{X\in Ob(J)}$ of 
$\ca C$-morphisms, such that for each $J$-morphism $f:X\to Y$ we have $\phi_X=F(f)\circ\phi_Y$. The co-cone 
$(M,\phi)$ is a {\it colimit} of the diagram $F$ if for every co-cone $(N,\psi)$ to $F$ there exists a 
unique $\ca C$-morphism $u:M\to N$ such that $\psi_X=\phi_X\circ u$ for all $X\in Ob(J)$. 
\enddefinition 

Again, it will be fruitful to review the case where $\ca C$ is the category of sets. 
Thus, let $F:J\to\ca C$ (with $\ca C=Sets$) be a $J$-shaped diagram. Let $\wh M$ be the disjoint union of 
the sets $F(X)$ for $X\in Ob(J)$. Let $\sim$ be 
the relation on $\wh M$ given by $a\sim b$ if there exists a $J$-morphism $f:X\to Y$ such that 
$a\in F(X)$, $b\in F(Y)$, and $b=(a)F(f)$ is the image of $a$ under $F(f)$. Let $\approx$ be the 
symmetrization of $\sim$ (so that $a\approx b$ if either $a\sim b$ or $b\sim a$). 
As $\sim$ is transitive, $\approx$ is then an equivalence relation.  
 Let $M$ be the set $\wh M/\approx$ of equivalence classes and 
let $\phi$ be the set of all $\phi_X:F(X)\to M$, where $\phi_X$ 
is the mapping which sends $a\in F(X)$ to the $\approx$-equivalence class of $a$ in $\wh M$.  
Then $(M,\phi)$ is a colimit of $F$. 

Next, take $\ca C$ to be the category $Set^*$ and let $F:J\to\ca C$ be a $J$-shaped diagram in $\ca C$. 
Here we take $\wh M$ to be the pointed set obtained from the disjoint union of the pointed sets $F(X)$ 
(over all $X\in Ob(J)$) by identifying base-points. For any $J$-morphism $f:X\to Y$, the 
morphism $F(f)$ of pointed sets sends the base-point of $F(X)$ to the base-point of $F(Y)$, and we may 
therefore define the equivalence relation $\approx$ on $\wh M$ as in the preceding paragraph. Take 
$M=\wh M/\approx$. Again, 
for each $X\in Ob(J)$ one has the mapping $\phi_X:X\to M$ which sends $a\in X$ to the 
$\approx$-equivalence class of $a$ in $\wh M$, and $(M,\phi)$ is a colimit of $F$.

In passing now to the case where $\ca C$ is the category of partial groups, we face the problem that, 
in general, there will be no partial normal subgroup, and no ``quotient" partial group corresponding to the 
equivalence relation $\approx$. For this reason, restrictions must be placed on the sort of diagrams  
$F:J\to Part$ that may be considered.

\proclaim {Theorem A.5} Let $Part$ be the category of partial groups, let $Set^*$ be the category of pointed 
sets, and let $J$ be a small category, and let $F:J\to Part$ be a $J$-shaped diagram. Assueme: 
\roster 

\item "{(1)}" For each  ordered pair $(X,Y)$ of objects of $J$, there exists at most one $J$-morphism 
$X\to Y$.  

\item "{(2)}" For each $J$-morphism $f:X\to Y$ the kernel of the homomorphism $F(f):F(X)\to F(Y)$ of 
partial groups is trivial. 

\endroster 
Let $F_0:J\to Sets^*$ be the composition of $F$ with the 
forgetful functor. Then there exists a colimit $(\ca M,\phi)$ of $F$, and the forgetful 
functor sends $(\ca M,\phi)$ to a limit of $F_0$. 
\endproclaim 

\demo {Proof} Let $\wh{\ca M}$ be the pointed set obtained as the disjoint union of all of the partial 
groups $F(X)$, for $X\in Ob(J)$, and with base-points identified. Define $\wh{\bold D}$ to be the 
disjoint union of the domains $\bold D(F(X))$. There is then a mapping $\wh\Pi:\wh{\bold D}\to\wh{\ca M}$ 
whose restriction to 
$\bold D(F(X))$ is the product $\Pi_X$ on $F(X)$. The union of the inversion maps on the partial groups 
$F(X)$ is an involutory bijection on $\wh{\ca M}$, and one may check that $\wh{\ca M}$ is a partial 
group via these structures. 

Let $\sim$ be the equivalence relation on $\wh{\ca M}$ given by $a\sim b$ if there exists a 
$J$-morphism $f:X\to Y$ such that $a\in F(X)$, $b\in F(Y)$, and $F(f):a\maps b$. Let $\approx$ be the 
equivalence relation given by symmetrizing $\sim$, and extend $\approx$ to an equivalence relation on 
$\bold W(\wh{\ca M})$ in the component-wise way. That is, if $u=(a_1,\cdots,a_m)$ and 
$v=(b_1,\cdots,b_n)$ are words in the alphabet $\wh{\ca M}$ then $u\approx v$ if and only if $m=n$ and 
$a_i\approx b_i$ for all $i$. For $a\in\wh{\ca M}$ we write $[a]$ for the $\approx$-class of $a$. Then 
the $\approx$-class of a word $u=(a_1,\cdots,a_n)$ is the word $([a_1],\cdots,[a_n])$.
Let $\ca M$ be the pointed set $\wh{\ca M}/\approx$ (whose base-point is the 
equivalence class of the base-point of $\wh{\ca M}$, and let $\bold D$ be the set of all words 
$([a_1],\cdots,[a_n])$ having a representative $(a_1,\cdots,a_n)\in\wh{\bold D}$.]  

Let $u=(a_1,\cdots,a_n)\in\wh{\bold D}$, and assume that there exists at least one index $k$ such that 
$a_k$ is not the identity element of $\wh{\ca M}$. Then there is a unique object $X$ of $J$ such that 
$u\in\bold D(F(X))$. Let also $v=(b_1,\cdots,b_n)\in\wh{\bold D}$, and assume $u\approx v$. Then (2) 
implies that $b_k$ is not the identity element of $\wh{\ca M}$, and there is a unique object 
$Y$ of $J$ with $v\in\bold D(F(Y))$. Let $i$ be any index from $1$ to $n$ such that not 
both $a_i$ and $b_i$ are identity elements. Then neither $a_i$ nor $b_i$ is an identity element, 
and (1) implies that either there is a unique $J$-morphism $f:X\to Y$ and $F(f):a_i\maps b_i$, or there 
is a unique $J$-morphism $g:Y\to X$ and $F(g):b_i\maps a_i$. If there exist both a $J$-morphism 
$f:X\to Y$ and a $J$-morphism $g:Y\to X$ then (1) implies that $f$ and $g$ are isomorphisms, and are 
inverse to each other, whence $F(f):a_i\maps b_i$ if and only if $F(g):b_i\maps a_i$. We may 
therefore assume without loss of generality that there exists a $J$-morphism $f:X\to Y$ and that 
$F(f)$ maps $u$ to $v$ component-wise. As $F(f)$ is a homomorphism $F(X)\to F(Y)$ of partial groups  
we then have $\wh\Pi(u)=\wh\Pi(v)$. We have thus shown that $\wh\Pi$ induces a mapping 
$\Pi:\bold D\to\ca M$, and the reader may check that $\Pi$ is a product, as defined in 1.1. 
If $\l:\ca H\to\ca K$ is a homomorphism of partial groups and $a\in\ca H$, then $(a\i)\l =(a\l)\i$, 
so there is a well-defined inversion mapping $\ca M\to\ca M$ given by $[a]\i=[a\i]$. Again it is left to 
the reader to check that with these structures, $\ca M$ is a partial group. 

For each $X\in Ob(J)$ define $\phi_X:F(X)\to\ca M$ by $a\phi_X=[a]$. Then $\phi_X$ is a homomorphism. If 
$a'\in F(X)$ with $[a]=[a']$ then $a=a'$, since the only $J$-morphism $X\to X$ is the identity morphism. 
Thus $\phi_X$ is injective. For any $J$-morphism $f:X\to Y$ and any $a\in F(X)$ we have $[(a)F(f)]=[a]$, 
so $(\ca M,\phi)$ is a co-cone of $F$. 

Let $(\ca N,\psi)$ be an arbitrary co-cone of $F$. Thus $\psi_X=F(f)\circ\psi_Y$ whenever 
$f:X\to Y$ is a $J$-morphism. That is, we have $a\psi_X=b\psi_Y$ if $F(f):a\maps b$, and thus 
there is a well-defined mapping $\s:\ca M\to\ca N$ given by $[a]\maps[a\psi_X]$ for $a\in F(X)$. 
Moreover, we have $\phi_X\circ\s=\psi_X$, and $\s$ is the unique such mapping $\ca M\to\ca N$. One 
checks that $\s$ is a homomorphism of partial groups, in order to complete the proof that 
$(\ca M,\phi)$ is a colimit of $F$. 
\qed 
\enddemo

\enddocument